\documentclass[oneside,english]{nlaauth}
\usepackage[T1]{fontenc}
\usepackage[latin9]{inputenc}
\usepackage{array}
\usepackage{multirow}
\usepackage{amsmath}
\usepackage{amssymb}
\usepackage{graphicx}
\usepackage{rotfloat}

\makeatletter

\providecommand{\tabularnewline}{\\}

\theoremstyle{remark}
\newtheorem{rem}{\protect\remarkname}

\def\volumeyear{2021}
\runningheads{Q. Chen, A. Ghai, and X. Jiao}{HILUCSI: Simple, Robust, and Fast Multilevel ILU}

\usepackage{times}
\usepackage[hidelinks]{hyperref}

\@ifundefined{showcaptionsetup}{}{%
 \PassOptionsToPackage{caption=false}{subfig}}
\usepackage{subfig}
\makeatother

\usepackage{babel}
\providecommand{\remarkname}{Remark}

\begin{document}
\title{HILUCSI: Simple, Robust, and Fast Multilevel ILU \\
for Large-Scale Saddle-Point Problems from PDEs}
\author{Qiao Chen,$^{1}$ Aditi Ghai,$^{1,2}$ Xiangmin Jiao$^{1}$\corrauth}
\address{$^{1}$Department of Applied Mathematics \& Statistics and Institute
for Advanced Computational Science, Stony Brook University, Stony
Brook, NY 11794, USA. \break\\
$^{2}$Current address: Cadence Design Systems Inc, San Jose, CA 95134,
USA.}
\corraddr{E-mail: xiangmin.jiao@stonybrook.edu.}
\begin{abstract}
Incomplete factorization is a widely used preconditioning technique
for Krylov subspace methods for solving large-scale sparse linear
systems. Its multilevel variants, such as ILUPACK, are more robust
for many symmetric or unsymmetric linear systems than the traditional,
single-level incomplete LU (or ILU) techniques. However, the previous
multilevel ILU techniques still lacked robustness and efficiency for
some large-scale saddle-point problems, which often arise from systems
of partial differential equations (PDEs). We introduce \emph{HILUCSI},
or \emph{Hierarchical Incomplete LU-Crout with Scalability-oriented
and Inverse-based dropping}. As a multilevel preconditioner, HILUCSI
statically and dynamically permutes individual rows and columns to
the next level for deferred factorization. Unlike ILUPACK, HILUCSI
applies symmetric preprocessing techniques at the top levels but always
uses unsymmetric preprocessing and unsymmetric factorization at the
coarser levels. The deferring combined with mixed preprocessing enabled
a unified treatment for nearly or partially symmetric systems and
simplified the implementation by avoiding mixed $1\times1$ and $2\times2$
pivots for symmetric indefinite systems. We show that this combination
improves robustness for indefinite systems without compromising efficiency.
Furthermore, to enable superior efficiency for large-scale systems
with millions or more unknowns, HILUCSI introduces a scalability-oriented
dropping in conjunction with a variant of inverse-based dropping.
We demonstrate the effectiveness of HILUCSI for dozens of benchmark
problems, including those from the mixed formulation of the Poisson
equation, Stokes equations, and Navier-Stokes equations. We also compare
its performance with ILUPACK and the supernodal ILUTP in SuperLU.
\end{abstract}
\maketitle

\section{Introduction}

Krylov subspace (KSP) methods, such as GMRES \cite{Saad03IMS,Saad86GMRES}
and BiCGSTAB \cite{vanderVorst92BiCGSTAB}, are widely used for solving
large-scale sparse unsymmetric or indefinite linear systems, especially
those arising from numerical discretizations of partial differential
equations (PDEs). For relatively ill-conditioned matrices, the KSP
methods can significantly benefit from a robust and efficient preconditioner.
Among these preconditioners, \emph{incomplete LU} (or\emph{ ILU})
is one of the most successful. Given a linear system $\boldsymbol{A}\boldsymbol{x}=\boldsymbol{b}$,
the ILU, or more precisely \emph{incomplete LDU} (or \emph{ILDU})
\emph{factorization }of $\boldsymbol{A}$ is an approximate factorization
\begin{equation}
\boldsymbol{P}^{T}\boldsymbol{A}\boldsymbol{Q}\approx\boldsymbol{L}\boldsymbol{D}\boldsymbol{U}.\label{eq:LU}
\end{equation}
On the right-hand side, $\boldsymbol{D}$ is a diagonal matrix, $\boldsymbol{L}$
is a unit lower triangular matrix, and $\boldsymbol{U}$ is an upper
triangular matrix; on the left-hand side, $\boldsymbol{P}$ and $\boldsymbol{Q}$
are row and column permutation matrices, respectively. Let $\boldsymbol{M}=\boldsymbol{L}\boldsymbol{D}\boldsymbol{U}$,
and $\boldsymbol{P}\boldsymbol{M}\boldsymbol{Q}^{T}$ is then a \emph{preconditioner}
of $\boldsymbol{A}$, or equivalently $\boldsymbol{M}$ is a preconditioner
of $\boldsymbol{P}^{T}\boldsymbol{A}\boldsymbol{Q}$. We consider
only right preconditioning in this work. Given the ILU factorization,
a right-preconditioned KSP method solves the preconditioned linear
system
\begin{equation}
\boldsymbol{A}\left(\boldsymbol{P}\boldsymbol{M}\boldsymbol{Q}^{T}\right)^{-1}\boldsymbol{y}=\boldsymbol{b},\label{eq:preconditioned-system}
\end{equation}
which ideally would converge much faster than solving the original
linear system, and then $\boldsymbol{x}=\left(\boldsymbol{P}\boldsymbol{M}\boldsymbol{Q}^{T}\right)^{-1}\boldsymbol{y}=\boldsymbol{Q}\boldsymbol{U}^{-1}\boldsymbol{D}^{-1}\boldsymbol{L}^{-1}\boldsymbol{P}^{T}\boldsymbol{y}$.

Earlier ILU methods, such as ILUTP \cite{chow1997experimental,Saad03IMS},
lack robustness for some indefinite systems (see, e.g., \cite{ernst2012difficult,lishao10,zhu2016generate}).
More recently, multilevel ILU (MLILU) techniques \cite{Boll06MPC,bank1999multilevel,saad2005multilevel},
a.k.a. \emph{multilevel block factorization} \cite{vassilevski2008multilevel},
have significantly improved the robustness of ILU for many applications.
A two-level preconditioner $\hat{\boldsymbol{M}}$ for the permuted
matrix $\boldsymbol{P}^{T}\boldsymbol{A}\boldsymbol{Q}$ can be constructed
via the approximation
\begin{equation}
\boldsymbol{P}^{T}\boldsymbol{A}\boldsymbol{Q}=\begin{bmatrix}\hat{\boldsymbol{B}} & \hat{\boldsymbol{F}}\\
\hat{\boldsymbol{E}} & \hat{\boldsymbol{C}}
\end{bmatrix}\approx\hat{\boldsymbol{M}}=\begin{bmatrix}\tilde{\boldsymbol{B}} & \tilde{\boldsymbol{F}}\\
\tilde{\boldsymbol{E}} & \hat{\boldsymbol{C}}
\end{bmatrix}=\begin{bmatrix}\boldsymbol{L}_{B} & \boldsymbol{0}\\
\boldsymbol{L}_{E} & \boldsymbol{I}
\end{bmatrix}\begin{bmatrix}\boldsymbol{D}_{B} & \boldsymbol{0}\\
\boldsymbol{0} & \boldsymbol{S}_{C}
\end{bmatrix}\begin{bmatrix}\boldsymbol{U}_{B} & \boldsymbol{U}_{F}\\
\boldsymbol{0} & \boldsymbol{I}
\end{bmatrix},\label{eq:mlilu-preconditioner}
\end{equation}
where $\hat{\boldsymbol{B}}\approx\tilde{\boldsymbol{B}}=\boldsymbol{L}_{B}\boldsymbol{D}_{B}\boldsymbol{U}_{B}$,
$\hat{\boldsymbol{E}}\approx\tilde{\boldsymbol{E}}=\boldsymbol{L}_{E}\boldsymbol{D}_{B}\boldsymbol{U}_{B}$,
and $\hat{\boldsymbol{F}}\approx\tilde{\boldsymbol{F}}=\boldsymbol{L}_{B}\boldsymbol{D}_{B}\boldsymbol{U}_{F}$.
The permutation matrices $\boldsymbol{P}$ and $\boldsymbol{Q}$ may
be obtained from some static reordering, static or dynamic pivoting,
or a combination of them. The Schur complement, i.e.,
\begin{equation}
\boldsymbol{S}_{C}=\hat{\boldsymbol{C}}-\tilde{\boldsymbol{E}}\tilde{\boldsymbol{B}}^{-1}\tilde{\boldsymbol{F}}=\hat{\boldsymbol{C}}-\boldsymbol{L}_{E}\boldsymbol{D}_{B}\boldsymbol{U}_{F},\label{eq:Schur-complement}
\end{equation}
is further approximately factorized recursively, resulting in a multilevel
ILU. In \cite{ghai2017comparison}, ILUPACK \cite{bollhofer2011ilupack,Boll06MPC},
which is a state-of-the-art MLILU package, was shown to be significantly
more robust than ILUTP and algebraic multigrid preconditioners \cite{hypre-user}
for indefinite systems. Nevertheless, the robustness and efficiency
of MLILU remained a challenge for saddle-point systems with millions
or more unknowns, which often arise from discretizations of Stokes,
Navier-Stokes, and Helmholtz equations, in computational fluid dynamics,
climate modeling, multiphysics coupling, etc.

The objective of this work is to improve the robustness and efficiency
of MLILU for saddle-point systems from PDEs. Our preconditioner was
motivated by two key observations. First, we observe that many linear
systems from systems of PDEs are \emph{``}nearly'' or ``partially''
symmetric, with some block structures. Without loss of generality,
assume matrix $\boldsymbol{A}\in\mathbb{R}^{n\times n}$ has the form
\begin{equation}
\boldsymbol{A}=\begin{bmatrix}\boldsymbol{B} & \boldsymbol{F}\\
\boldsymbol{E} & \boldsymbol{C}
\end{bmatrix}.\label{eq:block-matrix-form}
\end{equation}
It is worth noting that $\hat{\boldsymbol{B}}$ in (\ref{eq:mlilu-preconditioner})
may have different sizes from $\boldsymbol{B}$ in (\ref{eq:block-matrix-form})
due to permutation. For linear systems from PDEs, the nonzero pattern
of $\boldsymbol{A}$ is often \emph{nearly symmetric}, because in
some commonly used numerical methods (such as finite differences \cite{LeVeque07FDM},
finite elements \cite{Ciarlet2002FEM,ern2013theory}, or finite volumes
\cite{leveque2002finite}), the local support of the basis functions
(a.k.a. trial functions in finite elements) and that of the test functions
in the variational formulations are often the same or have significant
overlap. In addition, the numerical values are often \emph{nearly
symmetric} (i.e., $\left\Vert \boldsymbol{A}-\boldsymbol{A}^{T}\right\Vert \ll\left\Vert \boldsymbol{A}\right\Vert $),
because the numerical asymmetry is often due to small non-self-adjoint
terms (such as advection in a diffusion-dominant advection-diffusion
problem \cite[p. 243]{ern2013theory}) or due to truncation errors
(such as in a Petrov-Galerkin method for a self-adjoint PDE \cite[p. 88]{ern2013theory}).
For systems of PDEs, $\boldsymbol{A}$ may be \emph{partially symmetric}
in that $\left\Vert \boldsymbol{B}-\boldsymbol{B}^{T}\right\Vert \ll\left\Vert \boldsymbol{B}\right\Vert $
for the reasons mentioned above, but $\boldsymbol{E}$ and $\boldsymbol{F}^{T}$
differ significantly. Such partial symmetry may be due to strongly
imposed constraints in a variational formulation \cite{elman2014finite},
high-order treatment of Neumann boundary conditions in finite elements
\cite{cheung2019optimally} or finite differences \cite{LeVeque07FDM},
imposition of jump conditions in immersed/embedded boundary methods
\cite{johansen1998cartesian,peskin2002immersed}, etc. Second, we
observe that the state-of-the-art direct solvers, such as MUMPS \cite{amestoy2000mumps}
and PARDISO \cite{schenk2018PARDISO}, are highly optimized in terms
of cache performance, but they tend to scale superlinearly as the
problem sizes increase. In contrast, MLILU has poor locality due to
the dynamic nature of droppings, but its multilevel structure offers
additional opportunities to achieve near-linear time complexity while
being as robust as possible for very large systems, such as those
with millions or more unknowns.

Based on the preceding observations, we introduce a new preconditioner,
called \emph{HILUCSI }(pronounced as Hi-Luxi), which stands for \emph{Hierarchical
Incomplete LU-Crout with Scalability-oriented and Inverse-based dropping}.
As the name suggests, HILUCSI is a multilevel-ILU preconditioner that
utilizes the Crout version of ILU \cite{li2003crout}. In this aspect,
HILUCSI is closely related to ILUPACK \cite{bollhofer2011ilupack}.
The inverse-based dropping \cite{bollhofer2001robust,bollhofer2002relations,Boll06MPC}
in HILUCSI is also based on that of ILUPACK. However, HILUCSI improves
robustness and efficiency through a novel combination of several techniques.
First, HILUCSI is designed to take advantage of the near or partial
symmetry of the linear systems in a multilevel fashion. Specifically,
we apply symmetric preprocessing at the top levels for nearly or partially
symmetric matrices and apply unsymmetric factorization at lower levels
for all indefinite systems. This combination differs from other earlier
techniques for taking advantage of partial symmetry, such as using
$(\boldsymbol{A}+\boldsymbol{A}^{T})/2$ as an algebraic preconditioner
of $\boldsymbol{A}$ \cite{concus1976generalized,widlund1978lanczos}
or using the self-adjoint parts of the differential operators as a
physics-based preconditioner \cite{aksoylu2009family}.

Second, to construct its multilevel structure, HILUCSI introduces
a \emph{static} \emph{deferring} strategy to avoid nearly zero diagonal
entries, in conjunction with \emph{dynamic deferring} for avoiding
uncontrolled growth of the condition numbers of $\boldsymbol{L}_{B}$
and $\boldsymbol{U}_{B}$ in (\ref{eq:mlilu-preconditioner}) at each
level. Here, \emph{deferring} refers to a symmetric permutation operation
that delays some rows and their corresponding columns at one level
to the next level for deferred factorization. The dynamic deferring
in HILUCSI is similar to that in ILUPACK. Due to deferring, $\hat{\boldsymbol{B}}$
is in general smaller than $\boldsymbol{B}$ in HILUCSI. This behavior
makes HILUCSI different from, and potentially more robust than, a
straightforward block preconditioner for partially symmetric matrices
where $\hat{\boldsymbol{B}}$ is a simple permutation of $\boldsymbol{B}$
\cite{BenGolLie05NSS}. In addition, the static deferring in HILUCSI
simplifies its implementation for (nearly) symmetric saddle-point
problems, by eliminating the need of Bunch-Kaufman pivoting \cite{bunch1977some,Golub13MC}
as in \cite{greif2017sym,li2005crout} or similar $2\times2$ block
permutations as described in \cite{schenk2008large}.

Third, to achieve efficiency for large-scale problems, HILUCSI introduces
a \emph{scalability-oriented dropping}, which we use in conjunction
with inverse-based dropping. The primary goal of our scalability-oriented
dropping is to achieve (near) linear-time complexity in the number
of nonzeros in the input matrix. Although this goal shares some similarity
to the space-based droppings (such as those in ILUT \cite{saad1994ilut},
ICMF \cite{lin1999incomplete}, PARDISO \cite{schenk2018PARDISO})
as well as area-based dropping in \cite{lishao10}, it is different
in that linear-time complexity implies linear-space complexity but
not vice versa. We show that the scalability-oriented dropping along
with mixed symmetric and unsymmetric preprocessing enabled HILUCSI
to deliver superior robustness and efficiency compared to the previous
state of the art, such as ILUPACK \cite{bollhofer2011ilupack} and
supernodal ILUTP \cite{lishao10}, for saddle-point problems from
PDEs. 

The remainder of the paper is organized as follows. In Section~\ref{sec:background},
we review some background on incomplete LU factorization and its multilevel
variants. In Section~\ref{sec:HILUCSI}, we describe the algorithm
components of HILUCSI. Section~\ref{sec:Numerical-Results} presents
numerical results with HILUCSI as a right preconditioner for restarted
GMRES and compares its performance with some state-of-the-art packages.
Finally, Section~\ref{sec:Conclusions and Future Work} concludes
the paper with a discussion on future work. For completeness, we present
some analyses of stability and efficiency in the appendix.

\section{\label{sec:background}Preliminaries and related work}

In this section, we review some incomplete LU preconditioners. Because
there is a vast literature on preconditioning and ILU, we review only
some of the most relevant techniques and their mathematical analysis.
For comprehensive reviews, we refer readers to some surveys \cite{benzi2002preconditioning,chan1997approximate,saad2001iterative,wathen2015preconditioning}
and textbooks \cite{Saad03IMS,Van-der-Vorst:2003aa}.

\subsection{Single-level ILU}

The basic form of ILU in (\ref{eq:LU}), which we refer to as the
\emph{single-level ILU} (versus multilevel ILU), has been used to
solve linear systems from PDEs since the 1950s \cite{varga1959factorization}.
In 1977, Meijerink and van der Vorst \cite{meijerink1977iterative}
showed that incomplete Cholesky (IC) factorization is stable for a
symmetric M-matrix, i.e., a matrix with nonpositive off-diagonal entries
and nonnegative entries in its inverse. Since then, IC has been extended
to and become popular for SPD systems \cite{gupta2010adaptive,jones1995improved,lin1999incomplete,simon1988incomplete}.
However, ILU for unsymmetric or indefinite systems is significantly
more challenging, and it has been an active research topic over the
past few decades; see, e.g., \cite{chow1997experimental,lishao10,saad1994ilut,bollhofer2003robust,mayer2006alternative}.

\subsubsection{Variants of single-level ILU.}

In its simplest form, ILU does not involve any pivoting, and $\boldsymbol{L}$
and $\boldsymbol{U}$ preserve the sparsity patterns of the lower
and upper triangular parts of $\boldsymbol{P}^{T}\boldsymbol{A}\boldsymbol{Q}$,
respectively. This approach is often referred to as \emph{ILU0} or
\emph{ILU}(0). To improve its robustness, one may allow \emph{fills},
a.k.a. \emph{fill-ins}, which are new nonzeros in the $\boldsymbol{L}$
and $\boldsymbol{U}$ factors. The fills may be introduced based on
their levels in the elimination tree or based on the magnitude of
numerical values. The former leads to the so-called \emph{ILU}($k$),
which zeros out all the fills of level $k+1$ or higher in the elimination
tree. The combination of the two is known as \emph{ILU with dual thresholding}
(\emph{ILUT}) \cite{saad1994ilut}. The level-based fills may be replaced
with some other dropping to control the numbers of fills in each row
or column. The ILU algorithms in PETSc \cite{petsc-user-ref} and
hypre \cite{hypre-user} use some variants of ILUT.

ILUT may encounter zero or tiny pivots, which can lead to a breakdown
of the factorization. One may replace tiny pivots with a small value,
but such a trick is not robust \cite{chow1997experimental}. The robustness
may be improved by using pivoting, leading to the so-called \emph{ILUP}
\cite{saad1988preconditioning} and \emph{ILUTP} \cite{Saad03IMS}.
The\textit{ }ILU implementations in MATLAB \cite{MATLAB2019a}, SPARSKIT
\cite{Saad1994Sparsekit}, and the supernodal ILU in SuperLU\footnote{In this work, we use SuperLU to refer to its supernodal ILUTP, instead
of its better-known parallel complete LU factorization in SuperLU\_MT
or SuperLU\_Dist \cite{li2005overview}.} \cite{lishao10}, for example, are based on ILUTP. However, ILUTP
cannot prevent small pivots \cite{saad2005multilevel}, so it is still
not robust in practice; see, e.g., \cite{lishao10,mayer2007multilevel,saad2005multilevel}
and Section~\ref{sec:Numerical-Results} in this work for some failed
cases with ILUTP.

Another single-level ILU technique is the \emph{Crout-version of ILU},
a.k.a. \emph{ILUC} \cite{li2003crout,li2005crout}. In the context
of symmetric matrices, ILUC is also known as the \emph{left-looking
ILU} \cite{eisenstat1981algorithms}. At the $k$th step, ILUC updates
of the $k$th column in $\boldsymbol{L}$, which we denote by $\boldsymbol{\ell}_{k}$,
using the previous $k-1$ columns in $\boldsymbol{L}$, and it updates
the $k$th row of $\boldsymbol{U}$, which denoted by $\boldsymbol{u}_{k}^{T}$
, in a similar fashion. Unlike the standard $ijk$-ordered or ``right-looking''
ILU, ILUC updates $\boldsymbol{\ell}_{k}$ and $\boldsymbol{u}_{k}^{T}$
as late as possible. As a result, it allows dropping nonzeros in $\boldsymbol{\ell}_{k}$
and $\boldsymbol{u}_{k}^{T}$ as late as possible and avoids premature
dropping or partially updated columns. HILUCSI uses ILUC within each
level, along with some static and dynamic deferring described in Section~\ref{subsec:Static-and-dynamic}.

\subsubsection{\label{subsec:Accuracy-and-stability-SLILU}Accuracy and stability
of single-level preconditioner.}

For the ILU to be a ``robust'' preconditioner, it should be \emph{accurate}
and \emph{stable}. For the sake of simplicity, let us omit the permutation
matrices $\boldsymbol{P}$ and $\boldsymbol{Q}$ in the discussions,
and assume $\boldsymbol{M}=\boldsymbol{L}\boldsymbol{D}\boldsymbol{U}$
in (\ref{eq:LU}) is the preconditioner of $\boldsymbol{A}$. Following
\cite{benzi2002preconditioning}, we measure the accuracy and stability
of $\boldsymbol{M}$ by $\left\Vert \boldsymbol{A}-\boldsymbol{M}\right\Vert $
and $\left\Vert \boldsymbol{A}\boldsymbol{M}^{-1}-\boldsymbol{I}\right\Vert $,
respectively. In \cite{benzi2002preconditioning}, the Frobenius norm
was used for its ease of computation. We may use any other norm for
the convenience of analysis.

The accuracy of the preconditioner mostly depends on the dropping
strategies: Without dropping, $\left\Vert \boldsymbol{A}-\boldsymbol{M}\right\Vert $
is expected to be accurate up to machine precision. Most ILU techniques
utilize some form of \emph{dual thresholding} \cite{saad1994ilut},
by combining a numerical-value-based dropping and a combinatorial-structure-based
dropping. Several numerical-value-based dropping strategies have been
proposed for single-level ILU; see, e.g., \cite{bollhofer2001robust}
and \cite{mayer2006alternative}. These dropping strategies are often
used with a space-based dropping. To reduce space-based dropping,
some fill-reduction permutation, such as reverse Cuthill-Mckee (RCM)
\cite{chan1980linear} and approximate minimum degree (AMD) \cite{amestoy1996approximate},
are sometimes performed. In terms of stability, for $\left\Vert \boldsymbol{A}\boldsymbol{M}^{-1}-\boldsymbol{I}\right\Vert $
to be small (or bounded), $\Vert\boldsymbol{L}^{-1}\Vert$ and $\Vert\boldsymbol{U}^{-1}\Vert$
must be bounded by some small constant, as pointed out by Bollh\"ofer
\cite{bollhofer2001robust}. In addition, a tiny diagonal value in
$\boldsymbol{D}$ would also cause $\left\Vert \boldsymbol{A}\boldsymbol{M}^{-1}-\boldsymbol{I}\right\Vert $
to blow up, which corresponds to the danger of tiny pivots in ILU,
as emphasized by Chow and Saad in \cite{chow1997experimental}. To
improve the accuracy of ILU for some elliptic PDEs, Dupont et al.
\cite{dupont1968approximate} introduced \emph{modified ILU} (\emph{MILU}).
MILU modifies the diagonal entries to compensate for the discarded
entries so that the sum of each row is preserved. This technique appears
to be quite popular; see, e.g., \cite{gustafsson1978class,elman1986stability,lishao10}.
However, we do not use MILU, because it was ineffective for linear
systems from more general PDEs in our testing.

The dual requirements of accuracy and stability cause a dilemma for
single-level ILU for large-scale systems from PDEs, because the condition
number of $\boldsymbol{A}$, $\kappa(\boldsymbol{A})$, grows as the
number of unknowns $n$ increases. Specifically, it is well known
that for parabolic or elliptic PDEs, $\kappa(\boldsymbol{A})$ is
inversely proportional to $h^{-2}$ for some edge length measure $h$
\cite{LeVeque07FDM,ern2013theory}. Assuming isotropic meshes, we
expect $\kappa(\boldsymbol{A})$ to be proportional to $\mathcal{O}(n^{2/d})$,
where $d$ is the topological dimension of the domain. Assuming $\boldsymbol{A}$
is normalized such that $\Vert\boldsymbol{A}\Vert=\Theta(1)$, so
is $\Vert\boldsymbol{M}\Vert$. If $\boldsymbol{M}$ contains the
physically meaningful ``modes'' of the PDE that correspond to the
smallest singular values of $\boldsymbol{A}$, then $\left\Vert \boldsymbol{M}^{-1}\right\Vert \approx\Vert\boldsymbol{A}^{-1}\Vert$,
and $\left\Vert \boldsymbol{U}^{-1}\right\Vert \left\Vert \boldsymbol{D}^{-1}\right\Vert \left\Vert \boldsymbol{L}^{-1}\right\Vert \geq\left\Vert \boldsymbol{M}^{-1}\right\Vert \approx\mathcal{O}(n^{2/d})$.
Hence, it is challenging, if not impossible, to devise accurate single-level
ILU that is also stable (i.e., with bounded $\left\Vert \boldsymbol{L}^{-1}\right\Vert $,
$\left\Vert \boldsymbol{D}^{-1}\right\Vert $, and $\left\Vert \boldsymbol{U}^{-1}\right\Vert $)
for large-scale systems from PDEs. Note that although some preprocessing
techniques (such as equilibration \cite{VDS69CNE,golub2012matrix})
can be applied to reduce $\Vert\boldsymbol{A}^{-1}\Vert$ and $\left\Vert \boldsymbol{M}^{-1}\right\Vert $,
but they cannot alter this asymptotic behavior.  

\subsection{Multilevel ILU}

\subsubsection{Variants of MLILU.}

The aforementioned difficulties of single-level ILU are mitigated
by leveraging the multilevel structure in an MLILU. There are numerous
variants of MLILU in the literature, such as those in ILUM \cite{saad1996ilum},
BILUTM \cite{saad1999bilutm}, ARMS \cite{saad2002arms}, ILUPACK
\cite{bollhofer2011ilupack,Boll06MPC}, MDRILU \cite{zhang2001multilevel},
ILU++ \cite{mayer2007ilupp,mayer2007multilevel}, etc. Among these,
ILUPACK, MDRILU, and ILU++ aim to improve the robustness (and potentially
also the efficiency) of ILU using dynamic permutations. HILUCSI shares
similar goals as these preconditioners. Some methods use static reordering
to improve robustness in ILUT within a block, such as \cite{maclachlan2007greedy,osei2015matrix},
which we do not consider due to their lack of dynamic control of the
stability of triangular factors. Note that some MLILU techniques,
such as those in ILUM and BILUTM, were developed to expose parallelism,
and some others, such as those in \cite{bank1999incomplete,bank1999multilevel,bank1994hierarchical},
constructed an algebraic analogy of multigrid methods. The parallelization
of ILU and the mathematical connection between MLILU and multigrid
methods are beyond the scope of this work. 

\subsubsection{Accuracy and stability of MLILU preconditioner.}

Let us assume $\hat{\boldsymbol{M}}$ in (\ref{eq:mlilu-preconditioner})
is a preconditioner of $\boldsymbol{P}^{T}\boldsymbol{A}\boldsymbol{Q}$.
We need to compute $\hat{\boldsymbol{M}}^{-1}\boldsymbol{u}$ for
a block vector $\boldsymbol{u}=\begin{bmatrix}\boldsymbol{u}_{1}\\
\boldsymbol{u}_{2}
\end{bmatrix}$, i.e.,
\begin{equation}
\hat{\boldsymbol{M}}^{-1}\boldsymbol{u}=\begin{bmatrix}\tilde{\boldsymbol{B}}^{-1}\boldsymbol{u}_{1}\\
\boldsymbol{0}
\end{bmatrix}+\begin{bmatrix}-\tilde{\boldsymbol{B}}^{-1}\hat{\boldsymbol{F}}\\
\boldsymbol{I}
\end{bmatrix}\boldsymbol{S}_{C}^{-1}(\boldsymbol{u}_{2}-\hat{\boldsymbol{E}}\tilde{\boldsymbol{B}}^{-1}\boldsymbol{u}_{1}),\label{eq:Schur-complement-inverse}
\end{equation}
where $\boldsymbol{u}_{1}$ and $\boldsymbol{u}_{2}$ corresponds
to $\hat{\boldsymbol{B}}$ and $\hat{\boldsymbol{C}}$, respectively.
The stable computation of $\hat{\boldsymbol{M}}^{-1}\boldsymbol{u}$
requires the stable computation of $\tilde{\boldsymbol{B}}^{-1}\boldsymbol{u}_{1}$,
which in turn requires $\left\Vert \boldsymbol{L}_{B}^{-1}\right\Vert $,
$\left\Vert \boldsymbol{D}_{B}^{-1}\right\Vert $, and $\left\Vert \boldsymbol{U}_{B}^{-1}\right\Vert $
to be bounded by a constant at each level, analogous to the boundedness
of $\left\Vert \hat{\boldsymbol{B}}\tilde{\boldsymbol{B}}^{-1}-\boldsymbol{I}\right\Vert $
for single-level ILU. Unlike single-level ILU, MLILU can permute rows
and columns (such as permuting those in $\boldsymbol{B}$ that lead
to large condition numbers of $\boldsymbol{L}_{B}$ and $\boldsymbol{U}_{B}$
to the trailing block $\hat{\boldsymbol{C}}$) for deferred factorization,
as proposed in ILUPACK \cite{bollhofer2011ilupack,Boll06MPC}. HILUCSI
takes a similar approach to ensure the well-conditioning of $\tilde{\boldsymbol{B}}$,
except that it uses different permutation strategies than that of
ILUPACK for (nearly) symmetric indefinite systems. 

For the preconditioner $\hat{\boldsymbol{M}}$ to be stable, it is
clear that the computation of $\boldsymbol{S}_{C}^{-1}\boldsymbol{u}_{2}$
(or $\boldsymbol{S}_{C}^{-1}\boldsymbol{v}$ for $\boldsymbol{v}=\boldsymbol{u}_{2}-\hat{\boldsymbol{E}}\tilde{\boldsymbol{B}}^{-1}\boldsymbol{u}_{1}$)
should also be as stable as possible. To this end, we apply the MLILU
on $\boldsymbol{S}_{C}$ recursively until $\boldsymbol{S}_{C}$ is
small enough for a dense factorization. As in single-level ILU, the
stability of the computation can be improved by performing equilibration
\cite{VDS69CNE,golub2012matrix} before applying MLILU on $\boldsymbol{S}_{C}$.
Although these preprocessing strategies are used in virtually all
the MLILU techniques, HILUCSI differs from them in its combination
of symmetric preprocessing at the top levels with unsymmetric preprocessing
at the coarse levels.

In terms of the accuracy of MLILU, the concerns within each level
are similar to those of single-level ILU, as described in Section~\ref{subsec:Accuracy-and-stability-SLILU}.
For example, the dropping strategies in ILUPACK \cite{Boll06MPC}
and ILU++ \cite{mayer2007ilupp,mayer2007multilevel} are variants
of their single-level ILU in \cite{bollhofer2001robust} and \cite{mayer2006alternative},
respectively. The accuracy of the Schur complement in (\ref{eq:Schur-complement})
is an additional concern in MLILU. In \cite{Boll06MPC}, Bollh\"ofer
and Saad proposed to improve the accuracy by using a formulation due
to Tismenetsky \cite{tismenetsky1991new}, which, unfortunately, often
leads to excessive fills. A simpler alternative is to tighten the
dropping criteria for the Schur complement, as used in \cite{bollhofer2011ilupack,mayer2007multilevel,saad2002arms,zhang2001multilevel}.
In this work, we use the latter strategy.

\subsection{Near-linear time preconditioners}

In this work, we aim at devising an accurate and stable preconditioner
that has near-linear time complexity. In particular, the preconditioner
$\hat{\boldsymbol{M}}$ should be constructed in (near) linear time
in the number of nonzeros in the input matrix $\boldsymbol{A}$, and
$\boldsymbol{M}^{-1}\boldsymbol{u}$ can be computed in (near) linear
time for any $\boldsymbol{u}\in\mathbb{R}^{n}$. Note that ILU0 or
ILUTP with some space-based dropping may have linear space complexity,
but they may still have superlinear time complexity. Our objective
near-linear time complexity is similar to that of Hackbusch's hierarchical
matrices \cite{hackbusch2015hierarchical}. However, Hackbusch measures
the accuracy by some norm of $\boldsymbol{A}-\boldsymbol{M}$, without
taking into account the stability condition on $\boldsymbol{A}\boldsymbol{M}^{-1}-\boldsymbol{I}$
for preconditioners \cite{benzi2002preconditioning}. For linear systems
arising from finite difference or finite element discretizations for
PDEs, the number of nonzeros per row is typically bounded by a constant.
For such systems, our objective is to achieve near-linear time complexity
in the number of unknowns, similar to algebraic multigrid methods
for elliptic PDEs \cite{RS87AM}. This objective is also similar to
that of AMLI for SPD systems \cite{axelsson1990algebraic}. However,
we aim to solve unsymmetric and indefinite systems, for which the
state-of-the-art multigrid methods (such as BoomerAMG in hypre \cite{hypre-user})
are not robust \cite{ghai2017comparison}, and AMLI is inapplicable.

\section{\label{sec:HILUCSI}Hierarchical ILU with scalability-oriented droppings}

\begin{figure}
\centering{}\includegraphics[width=0.55\columnwidth]{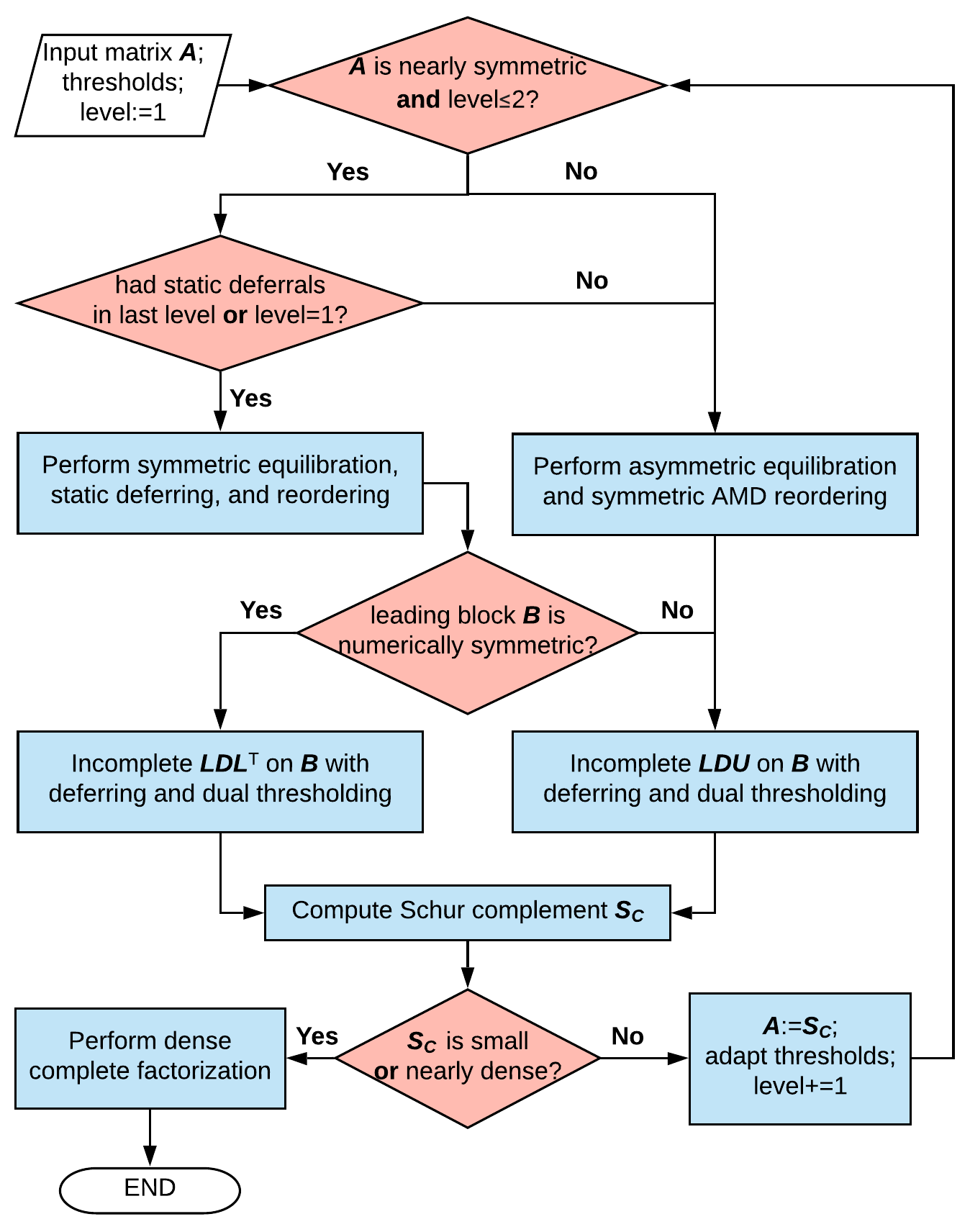}\caption{\label{fig:HILUCSI-factor}Flowchart of the multilevel ILDU factorization
in HILUCSI.}
\end{figure}

In this section, we describe the overall algorithm of HILUCSI. As
an MLILU, HILUCSI shares a similar control flow and some components
as others, such as ILUPACK and ILU++, except that we adapt the dropping
strategies to achieve better scalability for large-scale systems,
we adapt the pivoting/deferring strategies for simplicity for indefinite
systems without compromising robustness, and we adapt the preprocessing,
thresholds, and factorization techniques at different levels to take
into account near or partial symmetry of the systems. Figure~\ref{fig:HILUCSI-factor}
gives a schematic of the factorization procedure in HILUCSI, which
takes a sparse matrix $\boldsymbol{A}$, some thresholds for dropping
and permutations, along with a flag indicating near-symmetry as input.
In a nutshell, the algorithm dynamically builds a hierarchy of levels
by using static and dynamic deferring. Depending on whether the leading
block in the present level is symmetric or nearly symmetric, it performs
symmetric or unsymmetric preprocessing and factorization correspondingly,
and it performs a complete dense factorization at the coarsest level.
The algorithm returns the $\boldsymbol{L}_{B}$, $\boldsymbol{D}_{B}$,
$\boldsymbol{U}_{B}$, $\tilde{\boldsymbol{E}}$, and $\tilde{\boldsymbol{F}}$
of each level, which can then be used to compute (\ref{eq:Schur-complement-inverse}).
In the following, we focus on three key components of HILUCSI: the
scalability-oriented dropping with each level, the deferring strategies
for constructing the next level, and the mixed preprocessing strategies
for taking advantage of near and partial symmetry.

\subsection{\label{subsec:scalability-oriented-dropping}Scalability-oriented
dropping}

We first describe the core component of HILUCSI in the computation
at each level. For simplicity, let us omit the permutations (i.e.,
deferring) for now so that we can focus on the dropping strategy in
a two-level ILU by assuming $\boldsymbol{B}=\hat{\boldsymbol{B}}$.
Since one of the main objectives of HILUCSI is to achieve efficiency
for large-scale systems, we introduce a scalability-oriented dropping
specifically designed for MLILU. This dropping strategy has two parts.
First, we limit the number of nonzeros (nnz) in the $k$th column
of $\boldsymbol{L}$, namely $\boldsymbol{\ell}_{k}$, by a factor
of the nnz in the corresponding column in the input matrix, i.e.,
\begin{equation}
\text{nnz}(\boldsymbol{\ell}_{k})\leq\alpha\max\{\text{nnz}(\boldsymbol{a}_{k}),0.85\,\overline{\text{nnz}}(\boldsymbol{a}_{*})\},\label{eq:scalability-oriented-dropping}
\end{equation}
where $\alpha$ is a user-controllable parameter, and we refer to
it as the \emph{nnz factor}. The term $\overline{\text{nnz}}(\boldsymbol{a}_{*})$
denotes the average number of nonzeros in the columns of the input
matrix $\boldsymbol{A}$, and it is introduced to prevent excessive
dropping for columns with very few nonzeros in a highly nonuniform
matrix. Similarly, we limit the nnz in the rows of $\boldsymbol{U}$
in a similar fashion. Second, when computing the Schur complement
using (\ref{eq:Schur-complement}), we apply dropping to limit the
nnz in each column in $\boldsymbol{U}_{F}$ based on the right-hand
side of (\ref{eq:scalability-oriented-dropping}), and similarly for
each row in $\boldsymbol{L}_{E}$.
\begin{rem}
The two parts of the scalability-oriented dropping help control the
time complexity of updating $\boldsymbol{L}$ and $\boldsymbol{U}$
and computing the Schur complement, respectively. The first part is
easily incorporated into the Crout version of ILU. In particular,
at the $k$th step, after updating $\boldsymbol{\ell}_{k}$ (and $\boldsymbol{u}_{k}^{T}$),
we only keep up to $m$ largest-magnitude nonzeros, where $m$ is
equal to the right-hand side of (\ref{eq:scalability-oriented-dropping}).
The second part controls the time complexity of the sparse matrix-matrix
multiplication. It is important to note that in a multilevel setting,
$\text{nnz}(\boldsymbol{a}_{k})$ in (\ref{eq:scalability-oriented-dropping})
is the number of nonzeros in the $k$th column of the input matrix
$\boldsymbol{A}$ (assuming no permutation), instead of that of the
Schur complement from the previous level. For completeness, in Appendix~\ref{sec:Time-Complexity}
we present additional implementation details and a more detailed argument
of the linear time complexity of these two core steps within each
level.
\end{rem}
Besides the scalability-oriented dropping, HILUCSI also employs a
numerical-value-based strategy as a secondary dropping. We adopted
the inverse-based dropping as proposed in \cite{bollhofer2001robust}
and \cite{Boll06MPC}. More specifically, we estimate $\tilde{\kappa}_{L,k}\approx\left\Vert \boldsymbol{L}_{k}^{-1}\right\Vert _{\infty}$
and $\tilde{\kappa}_{U,k}\approx\left\Vert \boldsymbol{U}_{k}^{-1}\right\Vert _{1}$
incrementally as described in \cite{Higham87SCN}. Given a user-tunable
threshold $\kappa_{D}$ on $\left\Vert \boldsymbol{D}_{B}^{-1}\right\Vert _{\infty}$
for the level, we drop the $i$th nonzero in $\boldsymbol{\ell}_{k}$
if 
\begin{equation}
\kappa_{D}\tilde{\kappa}_{L,k}\vert\ell_{ik}\vert\leq\tau,\label{eq:inverse-based-thresholding_L}
\end{equation}
and we drop the $j$th nonzero in $\boldsymbol{u}_{k}^{T}$ if 
\begin{equation}
\kappa_{D}\tilde{\kappa}_{U,k}\vert u_{kj}\vert\leq\tau.\label{eq:inverse-based-thresholding-U}
\end{equation}
We refer to $\tau$ as the \emph{drop tolerance} (or in short, \emph{droptol}).
Like $\alpha$, $\tau$ is also a user-controllable threshold, and
it may vary from level to level. This inverse-based dropping can be
readily incorporated into the Crout version of ILU \cite{li2003crout}.
\begin{rem}
The inverse-based dropping described above was first proposed for
single-level ILU by Bollh\"ofer \cite{bollhofer2001robust}. It was
later adopted by Li et al. in Crout version of ILU \cite{li2003crout}
and then further adapted to multilevel ILU by Bollh\"ofer and Saad
in \cite{Boll06MPC}. Note that in \cite{Boll06MPC}, Bollh\"ofer
and Saad replaced $\tilde{\kappa}_{L,k}$ and $\tilde{\kappa}_{U,k}$
with a user-specified parameter $\kappa$, which is an upper-bound
of $\tilde{\kappa}_{L,k}$ and $\tilde{\kappa}_{U,k}$. Hence, our
inverse-based dropping is closer to its original form in \cite{bollhofer2001robust}
than that in \cite{Boll06MPC}. However, the scaling factor $\kappa_{D}$
in (\ref{eq:inverse-based-thresholding_L}) and (\ref{eq:inverse-based-thresholding-U})
was not present in \cite{bollhofer2001robust}. Although $\kappa_{D}$
may be blended into $\tau$, we find it conceptually clearer to separate
out $\kappa_{D}$, since $\kappa_{D}$ arises from the stability analysis
of the leading block as we summarize in Appendix~\ref{sec:Thresholding-in-inverse-based}.
$\kappa_{D}$ is also relevant in the deferring strategy for building
the multilevel structure, as we discuss next.
\end{rem}

\subsection{\label{subsec:Static-and-dynamic}Static and dynamic permutations
for deferred factorization}

The discussions in Section~\ref{subsec:scalability-oriented-dropping}
omitted permutations. It is well known that without permutations,
small pivots can make Gaussian elimination unstable \cite{Golub13MC},
and similar for ILU \cite{chow1997experimental}. However, pivoting
(such as partial pivoting in ILUTP \cite{saad1988preconditioning,Saad03IMS},
row pivoting in the supernodal ILUTP \cite{lishao10}, and dual pivoting
in ILU++ \cite{mayer2005ilucp}) requires sophisticated data structures
to implement. For symmetric indefinite systems, the Bunch-Kaufman
pivoting \cite{bunch1977some,Golub13MC} (such as in \cite{li2005crout}
and in \cite{schenk2008large}) incurs additional implementation complexities
by requiring permuting a combination of $1\times1$ and $2\times2$
pivots.

In HILUCSI, we exploit the multilevel structure to simplify both the
data structure and algorithm by only permuting the rows and the corresponding
columns in the leading block to the trailing blocks and deferring
them to the next level. We refer to this permutation strategy as \emph{deferring}.
In particular, we utilize two types of deferring. First, at the $k$th
step of Crout update, we dynamically permute the row and column to
the lower-right corner if the diagonal is smaller than a threshold
($\kappa_{D}$) or one of the estimated norms $\left\Vert \boldsymbol{L}_{B}^{-1}\right\Vert _{\infty}$
and $\left\Vert \boldsymbol{U}_{B}^{-1}\right\Vert _{1}$ exceeds
some threshold ($\kappa_{L}$ and $\kappa_{U}$). We refer to this
as \emph{dynamic deferring}, which is effective in resolving zero
or tiny pivots in most cases. In addition, we defer the zero and tiny
diagonal entries \emph{a priori}. We refer to this as \emph{static
pivoting}. In our experiments, we found that static pivoting is advantageous,
especially for saddle-point problems that have many zero diagonal
entries, probably because zero or tiny pivots tend to lead to rapid
growth of $\left\Vert \boldsymbol{L}_{B}^{-1}\right\Vert _{\infty}$
and $\left\Vert \boldsymbol{U}_{B}^{-1}\right\Vert _{1}$.

The static and dynamic deferring strategies in HILUCSI naturally result
in a multilevel structure, where within each level, we apply the algorithm
described in Section~\ref{subsec:scalability-oriented-dropping}.
To implement HILUCSI efficiently, we need a data structure that supports
efficient sequential access of the $k$th row of $\boldsymbol{L}$
along with all the rows in $\boldsymbol{U}_{1:k-1,k:n}$, and similarly
for the corresponding column in $\boldsymbol{U}$ and $\boldsymbol{L}_{k:n,1:k-1}$,
as required by the $k$th step of the Crout ILU. To this end, we use
a bi-index data structure based on that proposed by Li, Saad, and
Chow for the Crout version of ILU without pivoting for unsymmetric
matrices \cite{li2003crout}, which was based on that proposed by
Jones and Plassmann \cite{jones1995improved} for the row-version
of incomplete Cholesky factorization. The original data structures
in \cite{li2003crout} and \cite{jones1995improved} did not support
pivoting, but they require only about half of the storage and less
data movement than the more sophisticated data structures for ILU
with pivoting in \cite{li2005crout} and \cite{mayer2007multilevel}.
Since HILUCSI utilizes only static and dynamic deferring, we can extend
the more efficient data structure in \cite{chow1997experimental,li2003crout}
to support deferring without the extra memory overhead in \cite{li2005crout}
and \cite{mayer2007multilevel}. In particular, we allow the indices
to have a gap equal to the number of deferred rows and columns within
each level, and we eliminate the gap at the end of the ILU factorization
in the current level.

In terms of the thresholds, we use the same $\kappa_{D}$ for dynamic
deferring as that in (\ref{eq:inverse-based-thresholding_L}) and
(\ref{eq:inverse-based-thresholding-U}) for inverse-based dropping.
It typically suffices for $\kappa_{L}=\kappa_{U}=\kappa_{D}$, so
we can simply use $\kappa$ to denote their threshold and refer to
it as \emph{condest}. In terms of the thresholds, we observe from
numerical experimentation that it is desirable to tighten the thresholds
in level two by doubling $\alpha$, reducing $\tau$ by a factor of
10, and reducing $\kappa$ by a factor of two while restricting $\kappa\ge2$.
From level three, we revert $\alpha$ to the original value for efficiency,
while preserving the refined $\tau$ and $\kappa$. This choice is
because the dropping in $\boldsymbol{S}_{C}$ is amplified by $\tilde{\boldsymbol{B}}^{-1}$
in the term $-\tilde{\boldsymbol{B}}^{-1}\hat{\boldsymbol{F}}\boldsymbol{S}_{C}^{-1}\hat{\boldsymbol{E}}\tilde{\boldsymbol{B}}^{-1}$
in (\ref{eq:Schur-complement-inverse}), and the accuracy and stability
of level two appear to be the most important because it is the largest
Schur complement.
\begin{rem}
The dynamic deferring described above is similar to that for unsymmetric
matrices in \cite{Boll06MPC}. In \cite{schenk2008large}, Schenk
et al. described a dynamic deferring strategy for symmetric indefinite
systems, which permuted a combination of $1\times1$ and $2\times2$
blocks, motivated by the Bunch-Kaufman pivoting. The static deferring
in HILUCSI appears to be new. It eliminates the need of pivoting or
deferring of $2\times2$ blocks, and in turn, significantly simplifies
the treatment for symmetric indefinite systems without compromising
robustness or efficiency.
\end{rem}
It is worth noting that deferring is not foolproof. For example, if
the thresholds are too tight, then many rows and columns may be deferred,
leading to too many levels, undermining the robustness and also the
goal of near-linear time complexity. As another example, a matrix
may have all zero diagonals in the extreme case, and static deferring
would defer the whole matrix to the next level indefinitely. We mitigate
this issue by resorting to unsymmetric preprocessing and unsymmetric
factorization in coarse levels, as we address next.

\subsection{\label{subsec:Mixed-symmetric-and}Mixed symmetric and unsymmetric
factorization and preprocessing}

In HILUCSI, a novel feature is that it mixes symmetric and unsymmetric
techniques at different levels, to improve robustness and efficiency
for nearly or partially symmetric matrices. This mixed procedure has
two aspects. First, for matrices that are symmetric or partially symmetric,
i.e., the leading block $\boldsymbol{B}$ is symmetric in the input,
but $\boldsymbol{E}$ may or may not be equal to $\boldsymbol{F}^{T}$,
we apply symmetric factorization to $\boldsymbol{B}$, with static
and dynamic deferring. This combination allows HILUCSI to save the
factorization time by up to 50\% at the top levels for partially symmetric
matrices. However, we always resort to unsymmetric factorization at
the coarser levels. As alluded to in Section~\ref{subsec:Static-and-dynamic},
this mixed factorization allows HILUCSI to avoid indefinite deferring
when there are a large number of tiny diagonal entries in the input.
In our experiments, the benefits of improving robustness outweigh
the potential computational cost, since the coarser levels have small
sizes; see Section~\ref{subsec:Benefits-of-mixed} for more detail.

Secondly, HILUCSI employs different preprocessing techniques, including
reordering and equilibration, at different levels. As we noted in
Section~\ref{subsec:Accuracy-and-stability-SLILU}, fill-reduction
reordering (such as RCM and AMD) reduces the effect of droppings,
and equilibration improves the stability of the factorization by reducing
the condition numbers of the triangular factors. For nearly or partially
symmetric matrices, we apply reordering and equilibration symmetrically.
In particular, we apply RCM on $\text{nzp}(\boldsymbol{A})+\text{nzp}(\boldsymbol{A}^{T})$,
where $\text{nzp}$ denotes the nonzero pattern. In terms of equilibration,
we use MC64 \cite{duff1999design,duff2001algorithms} to compute the
row permutation vector $\boldsymbol{P}_{r}$ and the row and column
scaling vectors $\boldsymbol{D}_{r}$ and $\boldsymbol{D}_{c}$, respectively.
To preserve symmetry, we perform a post-processing step by setting
$\boldsymbol{P}_{r}=\boldsymbol{P}_{c}$ and $\tilde{\boldsymbol{D}}_{r}=\tilde{\boldsymbol{D}}_{c}=\sqrt{\boldsymbol{D}_{r}\boldsymbol{D}_{c}}$.
At coarser levels where we apply unsymmetric factorization, we always
use AMD reordering and MC64 equilibration directly. If the matrix
is unsymmetric, we apply reordering on $\text{nzp}(\boldsymbol{A})+\text{nzp}(\boldsymbol{A}^{T})$,
where denotes the nonzero pattern. In terms of the overall control
flow of preprocessing, we apply equilibration first, then static deferring,
and finally fill-reduction reordering on the leading block.
\begin{rem}
The mixed factorization and preprocessing improve the robustness for
nearly and partially symmetric matrices, especially those arising
from systems of PDEs, as we will show in Section~\ref{subsec:Benefits-of-mixed}.
We chose RCM and AMD for symmetric and unsymmetric reordering, respectively,
because it has previously been shown that RCM works better than AMD
for single-level ILU for symmetric matrices \cite{benzi1999orderings,gupta2010adaptive}.
We also observed similar behavior for multilevel ILU in our experiments,
probably because RCM tends to lead to smaller off-diagonal blocks,
which tends to improve the quality of the multilevel ILU. In terms
of symmetric equilibration, we use our own implementation of MC64
and the symmetrization process similar to that in HSL\_MC64 \cite{HSL_MC64}.
In \cite{duff2005strategies}, Duff and Pralet described a sophisticated
algorithm for symmetric indefinite systems, which involves $2\times2$
pivots and is more difficult to implement.
\end{rem}

\subsection{\label{subsec:Overall-time-complexity}Overall time complexity of
HILUCSI}

In terms of the computational cost, the core components of HILUCSI
achieve linear-time complexity at each level in the number of unknowns
for sparse linear systems from PDEs; see Appendix~\ref{sec:Time-Complexity}
for a detailed analysis. Note HILUCSI uses a complete dense factorization
in the last level, which has a cubic time complexity in its number
of rows and columns but excellent cache performance. To ensure the
overall linear-time complexity of factorization, we terminate the
recursion of HILUCSI when the final Schur complement is no more than
$\max\left\{ n^{1/3},C\right\} $ rows and columns, where $n$ is
the size of the original user input and $C$ is a constant. Base on
our experimentation, we found that $C=500$ leads to a negligible
dense-factorization cost. However, it is worth noting that the preprocessing
components in Section~\ref{subsec:Mixed-symmetric-and} may have
superlinear complexity in the worst case. In particular, AMD has quadratic-time
complexity in the worst case \cite{Heggernes01CCM,duff2001algorithms}.
In contrast, RCM has linear-time complexity \cite{chan1980linear},
which is yet another reason for using RCM instead of AMD at the top
levels. MC64 is superlinear in the worst case, but fortunately, it
has an expected linear-time complexity for most systems from PDEs
\cite{duff2001algorithms}. Hence, we claim that HILUCSI achieves
near-linear time complexity overall at each level. Note that the number
of levels in HILUCSI may grow as the problem size increases, albeit
very slowly. In Section~\ref{sec:Numerical-Results}, we will show
that HILUCSI indeed scales nearly linearly and performs better than
both ILUPACK and SuperLU for large systems with millions of unknowns.

\section{Numerical results\label{sec:Numerical-Results}}

We have implemented HILUCSI using the C++-11 standard. In this section,
we assess the robustness and efficiency of our implementation for
some challenging benchmark problems and compare its performance against
some state-of-the-art packages. In particular, we chose ILUPACK v2.4
\cite{bollhofer2011ilupack,ilupack} as the representative of multilevel
ILU, partially because HILUCSI is based on the same Crout-version
of multilevel ILU as in ILUPACK, and more importantly, ILUPACK has
been optimized for both unsymmetric and symmetric matrices. In comparison
with other packages, our tests showed that ILUPACK outperformed ARMS
in ITSOL v2 \cite{saad2002arms} by up to an order of magnitude for
larger unsymmetric systems. The improvement was even more significant
for symmetric systems. In our tests, ILUPACK was also significantly
more robust than ILU++ \cite{mayer2007ilupp,mayer2007multilevel}.
We chose the supernodal ILUTP in SuperLU v5.2.1 \cite{li2005overview,lishao10}
as a representative of the state-of-the-art single-level ILU. In all
of our tests, we used right-preconditioning for restarted GMRES, with
the dimension of the Krylov subspace limited to 30, i.e., GMRES(30).
We used $10^{-6}$ for the relative tolerance of GMRES and limited
the number of iterations to $500$. For HILUCSI and ILUPACK, we used
our implementation of flexible GMRES \cite{Saad03IMS}; for SuperLU,
we used GMRES implemented in PETSc v3.11.3. 

We conducted our tests on a single node of a cluster running CentOS
7.4 with two 2.5 GHz 12-core Intel Xeon E5-2680v3 processors and 64
GB of RAM. We compiled HILUCSI, SuperLU, and PETSc all using GCC 4.8.5
with the optimization option -O3, and we used the binary release of
ILUPACK v2.4 for GNU64. We accessed ILUPACK through its MATLAB mex
functions, of which the overhead is negligible. For accurate timing,
both turbo and power-saving modes were turned off for the processors.

\subsection{Baseline as ``black-box'' preconditioners}

As a baseline study, we assess HILUCSI for some benchmark problems
in the literature. We collected more than 60 larger-scale benchmark
problems that were highlighted in some recent ILU publications, which
were mostly from the SuiteSparse Matrix Collections \cite{davis2011university}
and the Matrix Market \cite{boisvert1997matrix}. For ill-conditioned
linear systems, we only consider those with a meaningful right-hand
side. We present results on some of the most challenging benchmark
problems that were highlighted in \cite{Boll06MPC}, \cite{lishao10},
and \cite{zhu2016generate}, together with two larger unsymmetric
systems for Navier-Stokes (N-S) equations. Table~\ref{tab:Benchmark-robustness}
summarizes these unsymmetric matrices, including their application
areas, types, sizes, and estimated condition numbers. Among the problems
omitted here, HILUCSI failed only for the system invextr1\_new in
\cite{zhu2016generate}, which has a large null-space dimension of
2,910 and also caused failures for all the methods tested in \cite{zhu2016generate}.
In addition, we generated two sets of symmetric indefinite systems
using FEniCS v2017.1.0 \cite{alnaes2015fenics} by discretizing the
3D Stokes equation and the mixed formulation of the Poisson equation.
These equations have a wide range of applications in computational
fluid dynamics (CFD), solid mechanics, heat transfer, etc. We discretized
the former using Taylor--Hood elements \cite{taylor1973numerical}
and discretized the latter using a mixture of linear Brezzi-Douglas-Marini
(BDM) elements \cite{brezzi1985two} and degree-0 discontinuous Galerkin
elements \cite{Cockburn2000}.  These problems are challenging because
the matrices have some nonuniform block structures, and they have
many zeros in the diagonals. To facilitate the scalability study,
for each set, we generated three linear systems using meshes of different
resolutions. Note that the matrices generated by FEniCS do not enforce
symmetry exactly and contain some nearly zero values due to rounding
errors. We manually filtered out the small values that are close to
machine precision and then enforced symmetry using $(\boldsymbol{A}+\boldsymbol{A}^{T})/2$.
 For this baseline comparison, we used droptol $\tau=10^{-4}$ for
all the codes, as in \cite{lishao10}, and used the recommended defaults
for the other parameters for most problems. For ILUPACK, we used MC64
matching, AMD ordering, and condest (i.e., $\kappa$) 5. For SuperLU,
when using its default options, we could solve four problems. We doubled
its ``fill factor'' from 10 to 20, which allowed SuperLU to solve
another five problems. For HILUCSI, we used nnz factor $\alpha=10$
and condest $\kappa=3$ for all the cases.

In Table~\ref{tab:Benchmark-robustness}, we report the overall runtimes
(including both factorization and solve times), numbers of GMRES iterations,
and the nnz ratios (i.e., the ratios of the number of nonzeros in
the output versus that in the input matrix, also known as the fill
ratio \cite{lishao10}) for each code. The fastest runtime for each
case is highlighted in bold. HILUCSI had a 95\% success rate for these
problems with the default parameters, and it was the fastest for 65\%
of the cases. For twotone, which is not a PDE-based problem, we could
not solve it unless we enlarge $\alpha$ to $15$. We note that for
all of the test cases, the final Schur complements in HILUCSI had
fewer than 500 rows and columns. ILUPACK solved 80\% of cases and
it was the fastest for 10\% of the cases. Among the failed cases,
ILUPACK ran out of the main memory for RM07R. For the symmetric problems,
ILUPACK automatically detects symmetric matrices and then applies
ILDL$^{T}$ factorization with mixed $1\times1$ and $2\times2$ pivots
automatically. This optimization in ILUPACK benefited its timing for
those problems but hurt its robustness for the two larger systems
from Stokes equations, which we could solve only by explicitly forcing
ILUPACK to use unsymmetric ILU. ILUPACK was unable to solve PR02R,
regardless of how we tuned its parameters. SuperLU was the least robust
among the three: It solved only 45\% of cases\footnote{In \cite{lishao10}, supernodal ILUTP had a higher success rate with
GMRES(50) and unlimited fill factor. We used GMRES(30) (the default
in PETSc \cite{petsc-user-ref}) and a fill factor 10 (the default
in SuperLU) or 20.}, and it was the fastest for 25\% of cases. Note that for the largest
system solved by all the codes, namely atmosmodl, HILUCSI outperformed
ILUPACK and SuperLU by a factor of 6 and 9, respectively. On the other
hand, for a medium-sized problem, namely e40r5000, SuperLU outperformed
HILUCSI and ILUPACK by a factor of 7.5 and 15, respectively. This
result shows that supernodal ILUTP excels in cache performance, but
its ILUTP is fragile compared to multilevel ILU in ILUPACK and HILUCSI.
Overall, HILUCSI delivered the best robustness and efficiency for
these cases.

\begin{sidewaystable}
\centering{}\caption{\label{tab:Benchmark-robustness}Baseline comparison of HILUCSI, ILUPACK,
and SuperLU, denoted as H, I, and S, respectively, using robust parameters
(droptol $10^{-4}$ for all and $\alpha=10$ for HILUCSI). $\times$,
and $-$ indicate failed factorization and stagnation of GMRES(30),
respectively. HILUCSI failed for twotone with $\alpha=10$, but it
took {\small{}7.71} seconds with $\alpha=15$. For SuperLU, `{*}'
indicates doubling fill-factor. The leaders are in bold.}
{\small{}}%
\begin{tabular}{c|c|c|c|c|c|c|c|c|c|c|c|c|c|c|c|c}
\hline 
\multirow{2}{*}{{\footnotesize{}Matrix}} & \multirow{2}{*}{{\footnotesize{}Appl.}} & \multirow{2}{*}{{\footnotesize{}n}} & \multirow{2}{*}{{\footnotesize{}nnz}} & \multirow{2}{*}{{\footnotesize{}Cond.}} & \multicolumn{3}{c|}{{\footnotesize{}factor. time (sec.)}} & \multicolumn{3}{c|}{{\footnotesize{}total time (sec.)}} & \multicolumn{3}{c|}{{\footnotesize{}GMRES iter.}} & \multicolumn{3}{c}{{\footnotesize{}nnz ratio}}\tabularnewline
\cline{6-17} \cline{7-17} \cline{8-17} \cline{9-17} \cline{10-17} \cline{11-17} \cline{12-17} \cline{13-17} \cline{14-17} \cline{15-17} \cline{16-17} \cline{17-17} 
 &  &  &  &  & {\footnotesize{}H} & {\footnotesize{}I} & {\footnotesize{}S} & {\footnotesize{}H} & {\footnotesize{}I} & {\footnotesize{}S} & {\footnotesize{}H} & {\footnotesize{}I} & {\footnotesize{}S} & {\footnotesize{}H} & {\footnotesize{}I} & {\footnotesize{}S}\tabularnewline
\hline 
\hline 
\multicolumn{17}{c}{{\footnotesize{}nearly symmetric, indefinite systems from PDEs}}\tabularnewline
\hline 
{\footnotesize{}venkat01} & {\footnotesize{}2D Euler} & {\footnotesize{}62,424} & {\footnotesize{}1,717,792} & {\footnotesize{}2.5e6} & {\footnotesize{}6.70} & {\footnotesize{}8.29} & \textbf{\footnotesize{}5.49} & {\footnotesize{}6.81} & {\footnotesize{}8.40} & \textbf{\footnotesize{}5.91} & {\footnotesize{}3} & {\footnotesize{}3} & {\footnotesize{}3} & {\footnotesize{}6.3} & {\footnotesize{}5.8} & {\footnotesize{}8.7}\tabularnewline
\hline 
{\footnotesize{}rma10} & \multirow{2}{*}{{\footnotesize{}3D CFD}} & {\footnotesize{}46,835} & {\footnotesize{}2,374,001} & {\footnotesize{}4.4e10} & {\footnotesize{}10.5} & {\footnotesize{}31.6} & \textbf{\footnotesize{}4.46} & {\footnotesize{}10.6} & {\footnotesize{}31.7} & \textbf{\footnotesize{}4.73} & {\footnotesize{}2} & {\footnotesize{}2} & {\footnotesize{}2} & {\footnotesize{}5.8} & {\footnotesize{}8.7} & {\footnotesize{}6.0}\tabularnewline
\cline{1-1} \cline{3-17} \cline{4-17} \cline{5-17} \cline{6-17} \cline{7-17} \cline{8-17} \cline{9-17} \cline{10-17} \cline{11-17} \cline{12-17} \cline{13-17} \cline{14-17} \cline{15-17} \cline{16-17} \cline{17-17} 
{\footnotesize{}mixtank\_new} &  & {\footnotesize{}29,957} & {\footnotesize{}1,995,041} & {\footnotesize{}2.2e11} & {\footnotesize{}39.7} & {\footnotesize{}128} & {\footnotesize{}$-$} & \textbf{\footnotesize{}40.2} & {\footnotesize{}128} & {\footnotesize{}$-$} & {\footnotesize{}9} & {\footnotesize{}3} & {\footnotesize{}$-$} & {\footnotesize{}15} & {\footnotesize{}41} & {\footnotesize{}$-$}\tabularnewline
\hline 
{\footnotesize{}Goodwin\_071} & {\footnotesize{}3D} & {\footnotesize{}56,021} & {\footnotesize{}1,797,934} & {\footnotesize{}1.4e7} & {\footnotesize{}20.2} & {\footnotesize{}15.4} & \textbf{\footnotesize{}6.45} & {\footnotesize{}20.8} & \textbf{\footnotesize{}15.5} & {\footnotesize{}16.3} & {\footnotesize{}13} & {\footnotesize{}3} & {\footnotesize{}81} & {\footnotesize{}15} & {\footnotesize{}8.2} & {\footnotesize{}9.2}\tabularnewline
\cline{1-1} \cline{3-17} \cline{4-17} \cline{5-17} \cline{6-17} \cline{7-17} \cline{8-17} \cline{9-17} \cline{10-17} \cline{11-17} \cline{12-17} \cline{13-17} \cline{14-17} \cline{15-17} \cline{16-17} \cline{17-17} 
{\footnotesize{}Goodwin\_095} & {\footnotesize{}Navier-} & {\footnotesize{}100,037} & {\footnotesize{}3,226,066} & {\footnotesize{}3.4e7} & {\footnotesize{}38.2} & {\footnotesize{}42.3} & \textbf{\footnotesize{}19.8} & {\footnotesize{}39.9} & {\footnotesize{}42.6} & \textbf{\footnotesize{}21.1}{\footnotesize{}{*}} & {\footnotesize{}18} & {\footnotesize{}3} & {\footnotesize{}4} & {\footnotesize{}16} & {\footnotesize{}9.2} & {\footnotesize{}12}\tabularnewline
\cline{1-1} \cline{3-17} \cline{4-17} \cline{5-17} \cline{6-17} \cline{7-17} \cline{8-17} \cline{9-17} \cline{10-17} \cline{11-17} \cline{12-17} \cline{13-17} \cline{14-17} \cline{15-17} \cline{16-17} \cline{17-17} 
{\footnotesize{}Goodwin\_127} & {\footnotesize{}Stokes} & {\footnotesize{}178,437} & {\footnotesize{}5,778,545} & {\footnotesize{}8.2e7} & {\footnotesize{}70.9} & {\footnotesize{}95.1} & \textbf{\footnotesize{}62.3} & {\footnotesize{}75.1} & {\footnotesize{}95.1} & \textbf{\footnotesize{}65.0}{\footnotesize{}{*}} & {\footnotesize{}24} & {\footnotesize{}3} & {\footnotesize{}4} & {\footnotesize{}16} & {\footnotesize{}12} & {\footnotesize{}15}\tabularnewline
\cline{1-1} \cline{3-17} \cline{4-17} \cline{5-17} \cline{6-17} \cline{7-17} \cline{8-17} \cline{9-17} \cline{10-17} \cline{11-17} \cline{12-17} \cline{13-17} \cline{14-17} \cline{15-17} \cline{16-17} \cline{17-17} 
{\footnotesize{}RM07R} & {\footnotesize{}(N-S)} & {\footnotesize{}381,689} & {\footnotesize{}37,464,962} & {\footnotesize{}2.2e11} & \textbf{\footnotesize{}3.1e3} & {\footnotesize{}$\times$} & {\footnotesize{}$\times$} & \textbf{\footnotesize{}3.3e3} & {\footnotesize{}$\times$} & {\footnotesize{}$\times$} & {\footnotesize{}75} & {\footnotesize{}$\times$} & {\footnotesize{}$\times$} & {\footnotesize{}44} & {\footnotesize{}$\times$} & {\footnotesize{}$\times$}\tabularnewline
\hline 
{\footnotesize{}PR02R} & \multirow{2}{*}{{\footnotesize{}2D N-S}} & {\footnotesize{}161,070} & {\footnotesize{}8,185,136} & {\footnotesize{}1.1e21} & \textbf{\footnotesize{}256} & {\footnotesize{}$-$} & {\footnotesize{}$-$} & \textbf{\footnotesize{}261} & {\footnotesize{}$-$} & {\footnotesize{}$-$} & {\footnotesize{}14} & {\footnotesize{}$-$} & {\footnotesize{}$-$} & {\footnotesize{}28} & {\footnotesize{}$-$} & {\footnotesize{}$-$}\tabularnewline
\cline{1-1} \cline{3-17} \cline{4-17} \cline{5-17} \cline{6-17} \cline{7-17} \cline{8-17} \cline{9-17} \cline{10-17} \cline{11-17} \cline{12-17} \cline{13-17} \cline{14-17} \cline{15-17} \cline{16-17} \cline{17-17} 
{\footnotesize{}e40r5000} &  & {\footnotesize{}17,281} & {\footnotesize{}553,956} & {\footnotesize{}2.5e10} & {\footnotesize{}18.0} & {\footnotesize{}36.8} & \textbf{\footnotesize{}2.26} & {\footnotesize{}18.7} & {\footnotesize{}36.8} & \textbf{\footnotesize{}2.4}{\footnotesize{}{*}} & {\footnotesize{}22} & {\footnotesize{}2} & {\footnotesize{}3} & {\footnotesize{}40} & {\footnotesize{}37} & {\footnotesize{}11}\tabularnewline
\hline 
{\footnotesize{}xenon2} & {\footnotesize{}materials} & {\footnotesize{}157,464} & {\footnotesize{}3,866,688} & {\footnotesize{}1.4e5} & \textbf{\footnotesize{}43.7} & {\footnotesize{}198} & {\footnotesize{}$-$} & \textbf{\footnotesize{}44.8} & {\footnotesize{}198} & {\footnotesize{}$-$} & {\footnotesize{}9} & {\footnotesize{}3} & {\footnotesize{}$-$} & {\footnotesize{}14} & {\footnotesize{}22} & {\footnotesize{}$-$}\tabularnewline
\hline 
{\footnotesize{}atmosmodl} & {\footnotesize{}Helmholtz} & {\footnotesize{}1,489,752} & {\footnotesize{}10,319,760} & {\footnotesize{}1.5e3} & \textbf{\footnotesize{}78.2} & {\footnotesize{}496} & {\footnotesize{}608} & \textbf{\footnotesize{}85.1} & {\footnotesize{}502} & {\footnotesize{}718} & {\footnotesize{}16} & {\footnotesize{}6} & {\footnotesize{}125} & {\footnotesize{}12} & {\footnotesize{}27} & {\footnotesize{}8.9}\tabularnewline
\hline 
\hline 
\multicolumn{17}{c}{{\footnotesize{}other unsymmetric systems}}\tabularnewline
\hline 
{\footnotesize{}onetone1} & {\footnotesize{}circuit} & {\footnotesize{}36,057} & {\footnotesize{}335,552} & {\footnotesize{}9.4e6} & \textbf{\footnotesize{}0.39} & {\footnotesize{}0.92} & {\footnotesize{}$-$} & \textbf{\footnotesize{}0.44} & {\footnotesize{}1.12} & {\footnotesize{}$-$} & {\footnotesize{}11} & {\footnotesize{}5} & {\footnotesize{}$-$} & {\footnotesize{}3.3} & {\footnotesize{}3.0} & {\footnotesize{}$-$}\tabularnewline
\cline{1-1} \cline{3-17} \cline{4-17} \cline{5-17} \cline{6-17} \cline{7-17} \cline{8-17} \cline{9-17} \cline{10-17} \cline{11-17} \cline{12-17} \cline{13-17} \cline{14-17} \cline{15-17} \cline{16-17} \cline{17-17} 
{\footnotesize{}twotone} & {\footnotesize{}simulation} & {\footnotesize{}120,750} & {\footnotesize{}1,224,224} & {\footnotesize{}4.5e9} & {\footnotesize{}$-$} & {\footnotesize{}18.3} & {\footnotesize{}$-$} & {\footnotesize{}$-$} & \textbf{\footnotesize{}18.5} & {\footnotesize{}$-$} & {\footnotesize{}$-$} & {\footnotesize{}2} & {\footnotesize{}$-$} & {\footnotesize{}$-$} & {\footnotesize{}12} & {\footnotesize{}$-$}\tabularnewline
\hline 
{\footnotesize{}bbmat} & {\footnotesize{}2D airfoil} & {\footnotesize{}38,744} & {\footnotesize{}1,771,722} & {\footnotesize{}5.4e8} & \textbf{\footnotesize{}31.4} & {\footnotesize{}36.4} & {\footnotesize{}57.2} & \textbf{\footnotesize{}31.9} & {\footnotesize{}36.5} & {\footnotesize{}59.0{*}} & {\footnotesize{}9} & {\footnotesize{}3} & {\footnotesize{}8} & {\footnotesize{}17} & {\footnotesize{}13} & {\footnotesize{}18}\tabularnewline
\hline 
\hline 
\multicolumn{17}{c}{{\footnotesize{}symmetric, saddle-point problems from PDEs}}\tabularnewline
\hline 
{\footnotesize{}S3D1} & \multirow{3}{*}{{\footnotesize{}3D Stokes}} & {\footnotesize{}18,037} & {\footnotesize{}434,673} & {\footnotesize{}1.2e7} & \textbf{\footnotesize{}1.63} & {\footnotesize{}6.40} & {\footnotesize{}7.68} & \textbf{\footnotesize{}1.63} & {\footnotesize{}6.56} & {\footnotesize{}9.90{*}} & {\footnotesize{}2} & {\footnotesize{}1} & {\footnotesize{}41} & {\footnotesize{}11} & {\footnotesize{}53} & {\footnotesize{}19}\tabularnewline
\cline{1-1} \cline{3-17} \cline{4-17} \cline{5-17} \cline{6-17} \cline{7-17} \cline{8-17} \cline{9-17} \cline{10-17} \cline{11-17} \cline{12-17} \cline{13-17} \cline{14-17} \cline{15-17} \cline{16-17} \cline{17-17} 
{\footnotesize{}S3D2} &  & {\footnotesize{}121,164} & {\footnotesize{}3,821,793} & {\footnotesize{}8.9e7} & \textbf{\footnotesize{}80.3} & {\footnotesize{}$\times$} & {\footnotesize{}$\times$} & \textbf{\footnotesize{}80.7} & {\footnotesize{}$\times$} & {\footnotesize{}$\times$} & {\footnotesize{}4} & {\footnotesize{}$\times$} & {\footnotesize{}$\times$} & {\footnotesize{}10} & {\footnotesize{}$\times$} & {\footnotesize{}$\times$}\tabularnewline
\cline{1-1} \cline{3-17} \cline{4-17} \cline{5-17} \cline{6-17} \cline{7-17} \cline{8-17} \cline{9-17} \cline{10-17} \cline{11-17} \cline{12-17} \cline{13-17} \cline{14-17} \cline{15-17} \cline{16-17} \cline{17-17} 
{\footnotesize{}S3D3} &  & {\footnotesize{}853,376} & {\footnotesize{}31,067,368} & {\footnotesize{}6.3e8} & \textbf{\footnotesize{}777} & {\footnotesize{}$\times$} & {\footnotesize{}$\times$} & \textbf{\footnotesize{}781} & {\footnotesize{}$\times$} & {\footnotesize{}$\times$} & {\footnotesize{}4} & {\footnotesize{}$\times$} & {\footnotesize{}$\times$} & {\footnotesize{}13} & {\footnotesize{}$\times$} & {\footnotesize{}$\times$}\tabularnewline
\hline 
{\footnotesize{}M3D1} & {\footnotesize{}3D} & {\footnotesize{}29,404} & {\footnotesize{}522,024} & {\footnotesize{}1.7e5} & \textbf{\footnotesize{}6.57} & {\footnotesize{}8.28} & {\footnotesize{}$-$} & \textbf{\footnotesize{}6.66} & {\footnotesize{}8.31} & {\footnotesize{}$-$} & {\footnotesize{}5} & {\footnotesize{}3} & {\footnotesize{}$-$} & {\footnotesize{}15} & {\footnotesize{}19} & {\footnotesize{}$-$}\tabularnewline
\cline{1-1} \cline{3-17} \cline{4-17} \cline{5-17} \cline{6-17} \cline{7-17} \cline{8-17} \cline{9-17} \cline{10-17} \cline{11-17} \cline{12-17} \cline{13-17} \cline{14-17} \cline{15-17} \cline{16-17} \cline{17-17} 
{\footnotesize{}M3D2} & {\footnotesize{}mixed} & {\footnotesize{}211,948} & {\footnotesize{}4,109,496} & {\footnotesize{}2.3e6} & \textbf{\footnotesize{}65.5} & {\footnotesize{}125} & {\footnotesize{}$\times$} & \textbf{\footnotesize{}67.1} & {\footnotesize{}127} & {\footnotesize{}$\times$} & {\footnotesize{}9} & {\footnotesize{}6} & {\footnotesize{}$\times$} & {\footnotesize{}16} & {\footnotesize{}28} & {\footnotesize{}$\times$}\tabularnewline
\cline{1-1} \cline{3-17} \cline{4-17} \cline{5-17} \cline{6-17} \cline{7-17} \cline{8-17} \cline{9-17} \cline{10-17} \cline{11-17} \cline{12-17} \cline{13-17} \cline{14-17} \cline{15-17} \cline{16-17} \cline{17-17} 
{\footnotesize{}M3D3} & {\footnotesize{}Poisson} & {\footnotesize{}1,565,908} & {\footnotesize{}31,826,184} & {\footnotesize{}3.8e7} & \textbf{\footnotesize{}570} & {\footnotesize{}1.2e3} & {\footnotesize{}$\times$} & \textbf{\footnotesize{}599} & {\footnotesize{}1.3e3} & {\footnotesize{}$\times$} & {\footnotesize{}21} & {\footnotesize{}12} & {\footnotesize{}$\times$} & {\footnotesize{}17} & {\footnotesize{}31} & {\footnotesize{}$\times$}\tabularnewline
\hline 
\end{tabular}
\end{sidewaystable}

\subsection{Optimized parameters for saddle-point problems}

The default parameters in the baseline comparison are robust for general
problems. However, they may be inefficient for saddle-point problems
from PDEs. We now compare the software for such problems, including
some nearly symmetric indefinite systems and purely symmetric saddle-point
problems.

\subsubsection{Nearly symmetric, indefinite systems.}

\begin{figure}
\subfloat[Relative speedup of factorization time.]{\includegraphics[width=0.45\columnwidth]{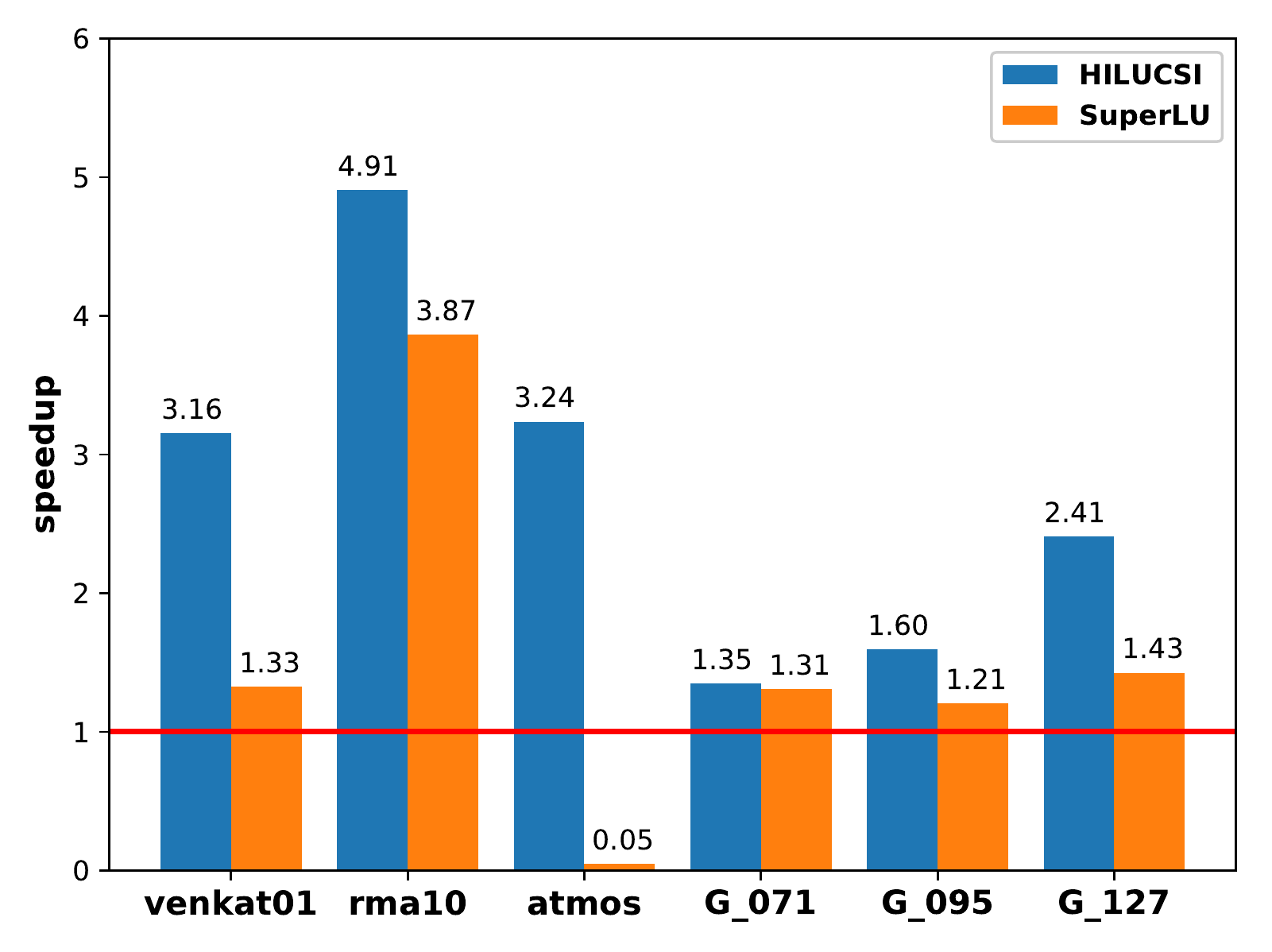}

}\hfill\subfloat[Relative speedup of total time.]{\includegraphics[width=0.45\columnwidth]{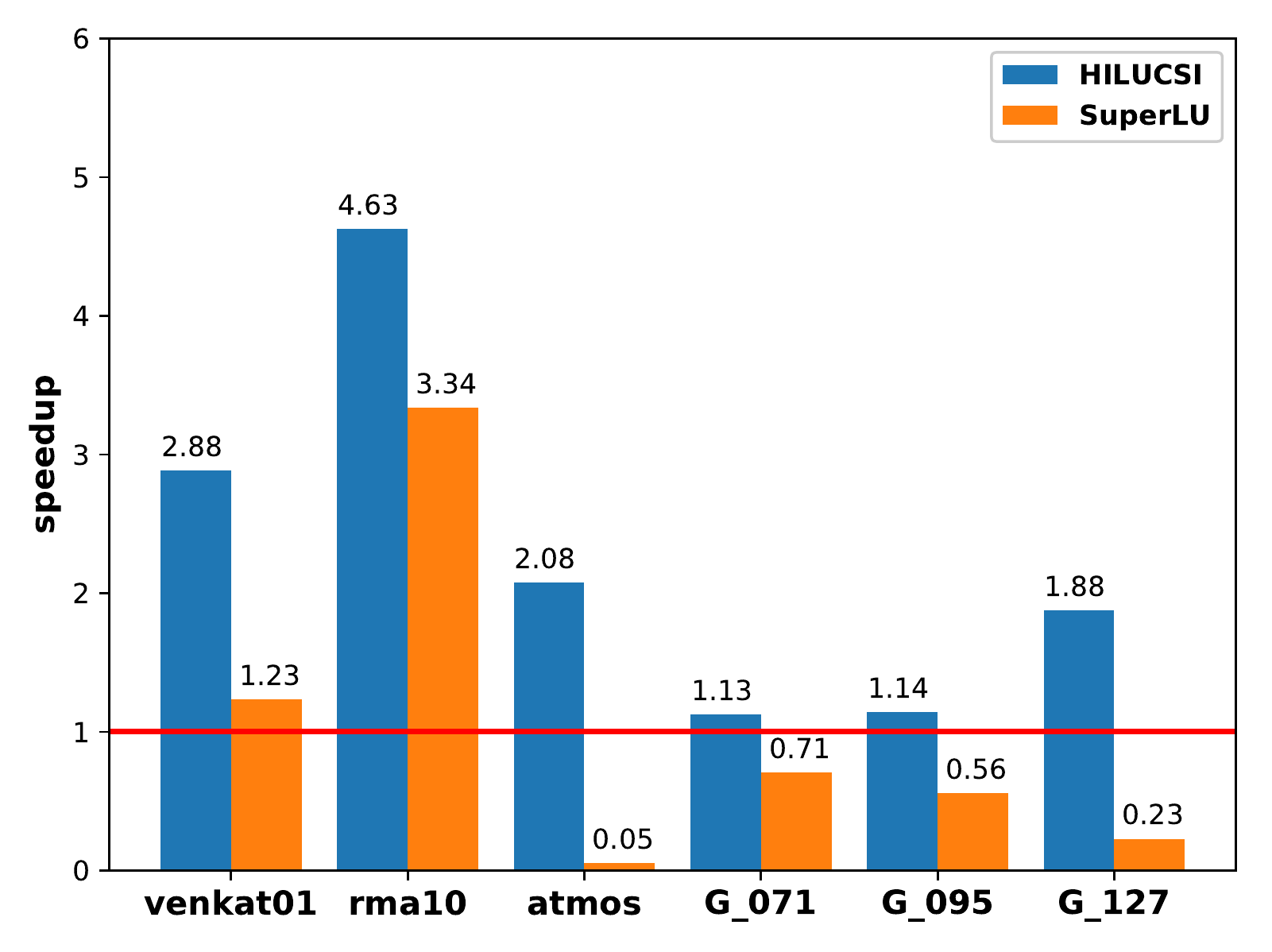}

}

\caption{\label{fig:speedups-asym-saddle-point}Speedups of (a) factorization
and (b) total times of HILUCSI and SuperLU versus ILUPACK for nearly
symmetric, indefinite problems with optimized parameters. Higher is
better. In the horizontal axis, ``atmos'' and ``G\_'' are short
for atmosmodl and Goodwin\_ in the horizontal axis, respectively.}
\end{figure}

For nearly symmetric matrices, we use six PDE-based problems in Table~\ref{tab:Benchmark-robustness},
which are from different types of equations in CFD, including 2D Euler,
3D Navier-Stokes equations, and Helmholtz equations. Table~\ref{tab:asym-saddle-problems}
shows the comparison of HILUCSI, ILUPACK, and SuperLU for these systems
in terms of the factorization times, total times, GMRES iterations,
and nonzero ratios. We highlighted the fastest runtimes in bold. For
a fair comparison, we used $\tau=0.01$ for all the codes, used $\kappa=5$
for both HILUCSI and ILUPACK, and used $\alpha=3$ for HILUCSI. It
can be seen that HILUCSI was the fastest for all the cases in terms
of both factorization and total times. Compared to ILUPACK, the lower
factorization cost of HILUCSI was due to a combination of smaller
fill factors, fewer levels, and lower time complexity (see Figure~\ref{fig:speedups-asym-saddle-point}).
However, HILUCSI required more GMRES iterations than ILUPACK, while
SuperLU required significantly more iterations for the largest systems.
In addition, we note that HILUCSI could solve all the problems with
$\alpha=2$, which improved the factorization time at the cost of
more GMRES iterations for some systems.

\begin{table}
\caption{\label{tab:asym-saddle-problems}Comparison of HILUCSI (denoted as
H) versus ILUPACK (I) and SuperLU (S) for nearly pattern-symmetric,
indefinite problems with optimized parameters (droptol $10^{-2}$
for all and $\alpha=3$ for HILUCSI). In matrix IDs, ``atmos'' and
``G\_'' are short for atmosmodl and Goodwin\_, respectively. Fastest
times are in bold.}

\centering{}{\scriptsize{}}%
\begin{tabular}{c|c|c|c|c|c|c|c|c|c|c|c|c|c|c}
\hline 
\multirow{2}{*}{{\scriptsize{}Matrix}} & \multicolumn{3}{c|}{{\scriptsize{}factor. time (sec.)}} & \multicolumn{3}{c|}{{\scriptsize{}total time (sec.)}} & \multicolumn{3}{c|}{{\scriptsize{}GMRES iters.}} & \multicolumn{3}{c|}{{\scriptsize{}nnz ratio}} & \multicolumn{2}{c}{{\scriptsize{}\#levels}}\tabularnewline
\cline{2-15} \cline{3-15} \cline{4-15} \cline{5-15} \cline{6-15} \cline{7-15} \cline{8-15} \cline{9-15} \cline{10-15} \cline{11-15} \cline{12-15} \cline{13-15} \cline{14-15} \cline{15-15} 
 & {\scriptsize{}H} & {\scriptsize{}I} & {\scriptsize{}S} & {\scriptsize{}H} & {\scriptsize{}I} & {\scriptsize{}S} & {\scriptsize{}H} & {\scriptsize{}I} & {\scriptsize{}S} & {\scriptsize{}H} & {\scriptsize{}I} & {\scriptsize{}S} & {\scriptsize{}H} & {\scriptsize{}I}\tabularnewline
\hline 
{\scriptsize{}venkat} & \textbf{\scriptsize{}1.11} & {\scriptsize{}3.50} & {\scriptsize{}2.64} & \textbf{\scriptsize{}1.25} & {\scriptsize{}3.62} & {\scriptsize{}2.94} & {\scriptsize{}7} & {\scriptsize{}5} & {\scriptsize{}7} & {\scriptsize{}3.0} & {\scriptsize{}2.8} & {\scriptsize{}2.4} & {\scriptsize{}3} & {\scriptsize{}3}\tabularnewline
\hline 
{\scriptsize{}rma10} & \textbf{\scriptsize{}2.31} & {\scriptsize{}11.3} & {\scriptsize{}2.93} & \textbf{\scriptsize{}2.47} & {\scriptsize{}11.4} & {\scriptsize{}3.43} & {\scriptsize{}9} & {\scriptsize{}4} & {\scriptsize{}9} & {\scriptsize{}2.0} & {\scriptsize{}3.8} & {\scriptsize{}2.8} & {\scriptsize{}5} & {\scriptsize{}6}\tabularnewline
\hline 
{\scriptsize{}atmos} & \textbf{\scriptsize{}10.5} & {\scriptsize{}33.8} & {\scriptsize{}686} & \textbf{\scriptsize{}19.8} & {\scriptsize{}41.0} & {\scriptsize{}748} & {\scriptsize{}33} & {\scriptsize{}22} & {\scriptsize{}75} & {\scriptsize{}2.9} & {\scriptsize{}4.0} & {\scriptsize{}6.2} & {\scriptsize{}3} & {\scriptsize{}3}\tabularnewline
\hline 
{\scriptsize{}G\_071} & \textbf{\scriptsize{}4.01} & {\scriptsize{}5.40} & {\scriptsize{}4.12} & \textbf{\scriptsize{}4.99} & {\scriptsize{}5.63} & {\scriptsize{}7.95} & {\scriptsize{}41} & {\scriptsize{}12} & {\scriptsize{}52} & {\scriptsize{}4.8} & {\scriptsize{}3.5} & {\scriptsize{}5.0} & {\scriptsize{}5} & {\scriptsize{}7}\tabularnewline
\hline 
{\scriptsize{}G\_095} & \textbf{\scriptsize{}7.36} & {\scriptsize{}11.8} & {\scriptsize{}9.74} & \textbf{\scriptsize{}10.7} & {\scriptsize{}12.3} & {\scriptsize{}21.9} & {\scriptsize{}78} & {\scriptsize{}14} & {\scriptsize{}84} & {\scriptsize{}4.8} & {\scriptsize{}3.9} & {\scriptsize{}5.7} & {\scriptsize{}6} & {\scriptsize{}7}\tabularnewline
\hline 
{\scriptsize{}G\_127} & \textbf{\scriptsize{}13.3} & {\scriptsize{}32.1} & {\scriptsize{}22.5} & \textbf{\scriptsize{}17.7} & {\scriptsize{}33.3} & {\scriptsize{}146} & {\scriptsize{}56} & {\scriptsize{}16} & {\scriptsize{}449} & {\scriptsize{}4.8} & {\scriptsize{}5.0} & {\scriptsize{}6.1} & {\scriptsize{}6} & {\scriptsize{}8}\tabularnewline
\hline 
\end{tabular}{\scriptsize\par}
\end{table}

Figure~\ref{fig:speedups-asym-saddle-point} shows the relative speedups
of HILUCSI and SuperLU versus ILUPACK in terms of factorization and
solve times. It can be seen that HILUCSI outperformed ILUPACK for
all six cases by a factor between 1.1 and 4.9. For the Goodwin problems,
it is clear that the relative speedup increased as the problem sizes
increased, thanks to the near-linear time complexity of HILUCSI as
discussed in Section~\ref{subsec:Overall-time-complexity}. We note
that ILUPACK has a parameter \textsf{elbow} for controlling the size
of reserved memory, but the parameter made no difference in our testing.
ILUPACK also has another parameter \textsf{lfil} for space-based dropping,
of which the use is discouraged in its documentation. Our tests showed
that using a small \textsf{lfil} in ILUPACK decreased its robustness,
while its time complexity was still higher than HILUCSI. 

We observe that although SuperLU outperformed ILUPACK in terms of
factorization times for all the Goodwin problems, it underperformed
in terms of the overall times for these problems, due to the slow
convergence of GMRES. This result again shows the superior robustness
of MLILU in HILUCSI and ILUPACK versus the single-level ILUTP in SuperLU.

\subsubsection{\label{subsec:Symmetric-saddle-point-problems}Symmetric saddle-point
problems.}

We now assess the robustness and efficiency of HILUCSI as the problem
sizes increase. To this end, we use the symmetric saddle-point problems
and compare HILUCSI with two different solvers in ILUPACK for symmetric
and unsymmetric matrices, respectively. Because supernodal ILUTP failed
for most of these problems, we do not include it in this comparison.
For these saddle-point problems, because there are static deferring,
our algorithm enabled symmetric matching in HILUCSI on the first two
levels, and we applied RCM for the first level and applied AMD ordering
for all the other levels. For ILUPACK, we used AMD ordering, as recommended
by ILUPACK's documentation.

Table~\ref{tab:sym-saddle-problems} shows the comparison of HILUCSI
with ILUPACK in terms of the numbers of GMRES iterations, the nonzero
ratios, and the numbers of levels, along with the runtimes of HILUCSI.
It is worth noting that symmetric ILUPACK failed for the two larger
systems for the Stokes equations due to encountering a structurally
singular system during preprocessing. For the two larger cases for
the mixed formulation of the Poisson equation, symmetric ILUPACK was
notably less robust and required many more GMRES iterations. Among
the four solved problems, symmetric ILUPACK improved the runtimes
of unsymmetric ILUPACK by a factor of 1.2 to 2.6, because the symmetric
version performed computations only on the lower triangular part and
used different dropping strategies. 
\begin{rem}
The results of unsymmetric versus symmetric ILUPACK in Table~\ref{tab:sym-saddle-problems}
show that it is sometimes more robust to solve symmetric indefinite
systems using unsymmetric solvers. Conventionally, it is believed
that ``it is rarely sensible to throw away symmetry in preconditioning''
\cite{wathen2015preconditioning}. Such conventional wisdom seems
to focus on the efficiency of single-level ILU, because it may reduce
the computational cost by up to half by using symmetric (versus unsymmetric)
factorization \cite{duff2005strategies,schenk2018PARDISO}. The use
of unsymmetric factorization at the coarse levels in HILUCSI is a
crucial reason for its robustness for symmetric saddle-point problems.
\end{rem}
\begin{table}
\caption{\label{tab:sym-saddle-problems}Comparison of HILUCSI (denoted as
H) with unsymmetric and symmetric ILUPACK (denoted by IU and IS, respectively)
for symmetric saddle-point systems with droptol $10^{-2}$ for all
and $\alpha=3$ for HILUCSI. $\times$ indicates failed factorization.}

\centering{}%
\begin{tabular}{c|c|c|c|c|c|c|c|c|c|c|c}
\hline 
\multirow{2}{*}{Matrix} & \multicolumn{2}{c|}{HILUCSI} & \multicolumn{3}{c|}{GMRES iters.} & \multicolumn{3}{c|}{nnz ratio} & \multicolumn{3}{c}{\#levels}\tabularnewline
\cline{2-12} \cline{3-12} \cline{4-12} \cline{5-12} \cline{6-12} \cline{7-12} \cline{8-12} \cline{9-12} \cline{10-12} \cline{11-12} \cline{12-12} 
 & factor. & total & H & IU & IS & H & IU & IS & H & IU & IS\tabularnewline
\hline 
S3D1 & 0.44 & 0.45 & 4 & 3 & 7 & 2.0 & 4.6 & 6.4 & 3 & 6 & 5\tabularnewline
\hline 
S3D2 & 5.56 & 5.83 & 7 & 3 & $\times$ & 2.5 & 6.3 & $\times$ & 3 & 6 & $\times$\tabularnewline
\hline 
S3D3 & 61.7 & 64.1 & 7 & 4 & $\times$ & 2.7 & 8.4 & $\times$ & 4 & 9 & $\times$\tabularnewline
\hline 
M3D1 & 0.69 & 0.78 & 14 & 6 & 11 & 2.7 & 7.9 & 5.8 & 4 & 8 & 5\tabularnewline
\hline 
M3D2 & 6.25 & 7.75 & 26 & 11 & 29 & 2.6 & 9.5 & 7.3 & 5 & 10 & 5\tabularnewline
\hline 
M3D3 & 52.9 & 76.8 & 53 & 24 & 62 & 2.6 & 10 & 7.2 & 6 & 15 & 5\tabularnewline
\hline 
\end{tabular}
\end{table}

\begin{figure}
\subfloat[Relative speedup of factorization time.]{\includegraphics[width=0.45\columnwidth]{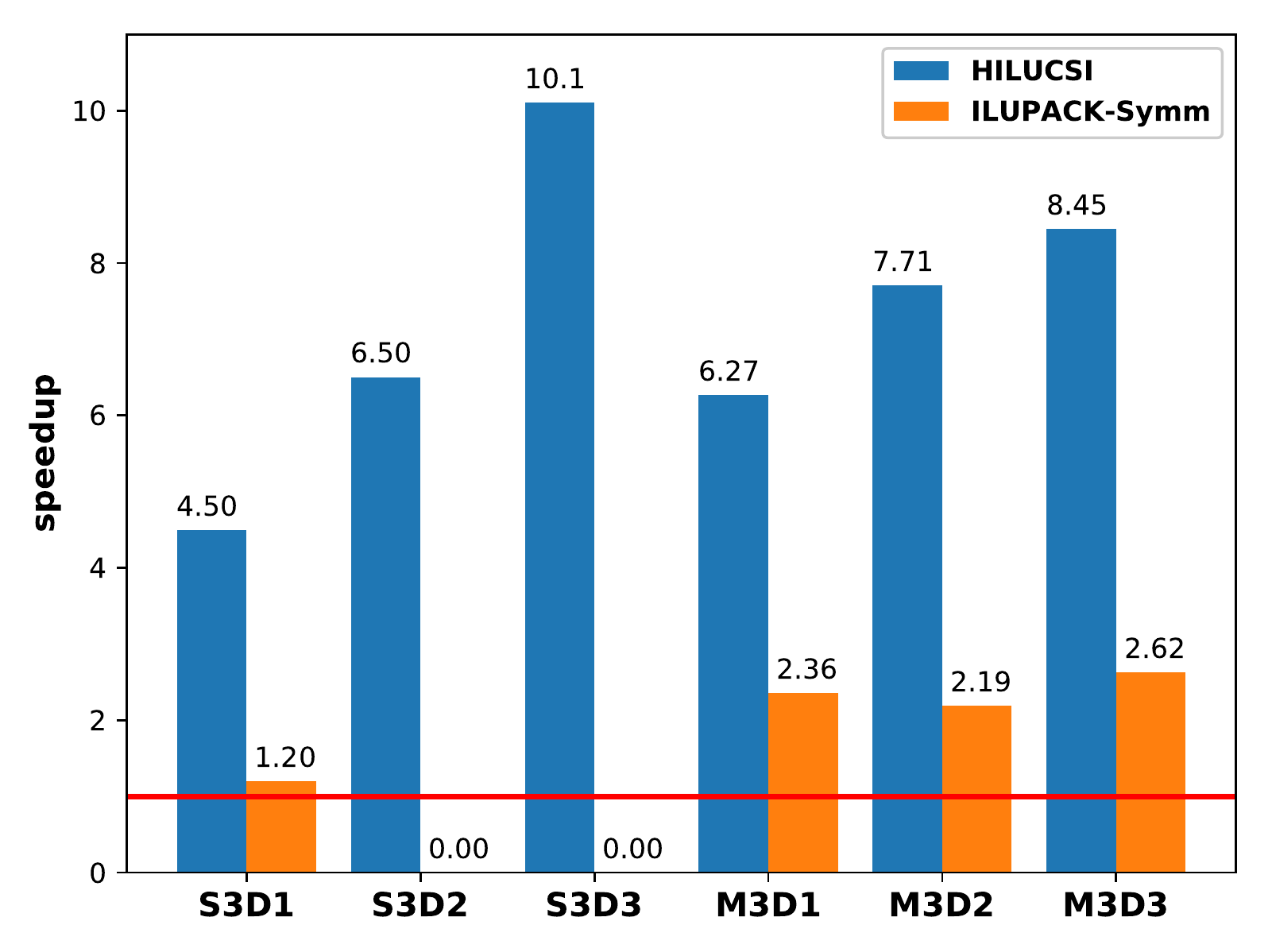}

}\hfill\subfloat[Relative speedup of total time.]{\includegraphics[width=0.45\columnwidth]{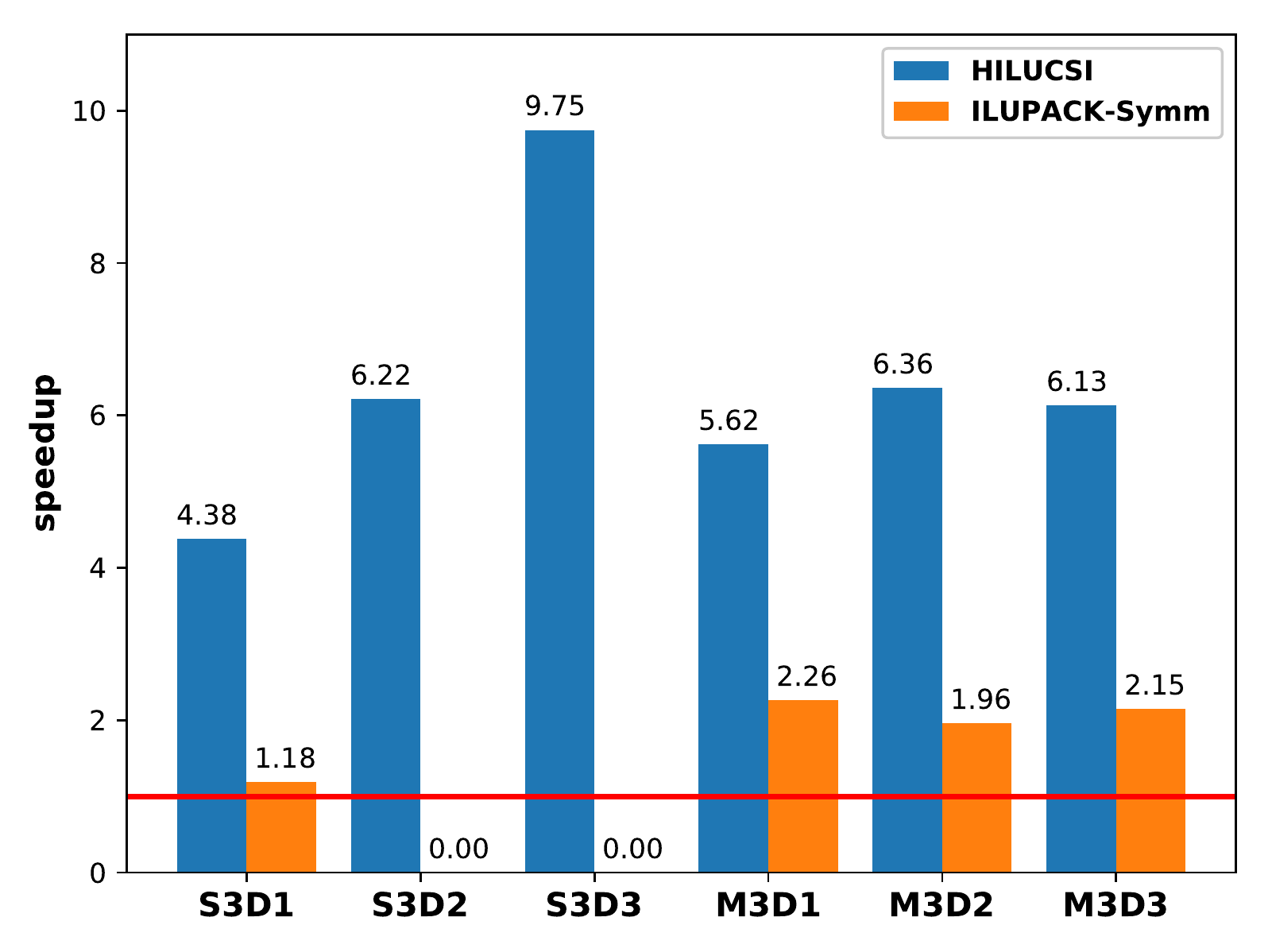}

}

\caption{\label{fig:speedups-saddle-points}Speedups of (a) factorization and
(b) total times of HILUCSI and symmetric factorization in ILUPACK
versus unsymmetric ILUPACK for symmetric saddle-point problems. Higher
is better.}
\end{figure}

In terms of efficiency, Figure~\ref{fig:speedups-saddle-points}
shows the relative speedups of HILUCSI and symmetric ILUPACK relative
to the unsymmetric ILUPACK. It can be seen that HILUCSI outperformed
the unsymmetric ILUPACK by a factor of four to ten for these problems.
The improvement was mostly due to the improved dropping in HILUCSI.
HILUCSI also had fewer levels than unsymmetric ILUPACK. Note that
the timing results in Table~\ref{tab:sym-saddle-problems} for HILUCSI
did not use symmetric factorization at any level. The use of symmetric
factorization in the first two levels further improved its overall
performance by 10--20\%. Note that in Figure~\ref{fig:speedups-saddle-points},
the relative speedup of HILUCSI versus ILUPACK grew as the problem
sizes increased. Hence, the better efficiency of HILUCSI for large-scale
systems is primarily due to its better scalability with respect to
problem sizes, thanks to its scalability-oriented dropping.

\subsection{\label{subsec:Benefits-of-mixed}Benefits of mixed preprocessing}

To assess the effectiveness of mixing symmetric and unsymmetric preprocessing
for HILUCSI as we described in Section~\ref{subsec:Mixed-symmetric-and},
we applied symmetric preprocessing on zero, one, and two levels. Table~\ref{tab:mixed-proprocessing}
shows a comparison of the factorization times, total times, GMRES
iterations, and nnz ratios for three different classes of problems.
It can be seen that for matrices with fully unsymmetric structures,
the use of symmetric preprocessing did not improve robustness and
even decreased efficiency. However, for many unsymmetric matrices
with nearly symmetric structures, using symmetric preprocessing on
the first level significantly improved robustness and efficiency.
The behavior is probably because the matrices from systems of PDEs
tend to have some block diagonal dominance as defined in \cite{feingold1962block},
which may be destroyed by unsymmetric permutations. Symmetric reordering
and equilibration permute the dominant block diagonal within a narrower
band, so that it may preserve block diagonal dominance more effectively.
Furthermore, when static deferring is invoked in (nearly) symmetric
saddle-point problems, using two levels of symmetric preprocessing
further reduced the factorization times, but the total runtime remained
about the same.

\begin{table}
\caption{\label{tab:mixed-proprocessing}Effect of mixing symmetric and unsymmetric
processing in HILUCSI. H0, H1, and H2 denote using zero, one, and
two levels of symmetric preprocessing.}

\centering{}{\small{}}%
\begin{tabular}{c|c|c|c|c|c|c|c|c|c|c|c|c}
\hline 
\multirow{2}{*}{{\small{}Matrix}} & \multicolumn{3}{c|}{{\small{}factor. time}} & \multicolumn{3}{c|}{{\small{}total time}} & \multicolumn{3}{c|}{{\small{}GMRES iters.}} & \multicolumn{3}{c}{{\small{}nnz ratio}}\tabularnewline
\cline{2-13} \cline{3-13} \cline{4-13} \cline{5-13} \cline{6-13} \cline{7-13} \cline{8-13} \cline{9-13} \cline{10-13} \cline{11-13} \cline{12-13} \cline{13-13} 
 & {\small{}H0} & {\small{}H1} & {\small{}H2} & {\small{}H0} & {\small{}H1} & {\small{}H2} & {\small{}H0} & {\small{}H1} & {\small{}H2} & {\small{}H0} & {\small{}H1} & {\small{}H2}\tabularnewline
\hline 
\hline 
\multicolumn{13}{c}{{\small{}general unsymmetric systems}}\tabularnewline
\hline 
{\small{}bbmat} & \textbf{\small{}31.4} & {\small{}45.5} & {\small{}55.2} & {\small{}31.9} & {\small{}46.3} & {\small{}55.9} & {\small{}9} & {\small{}11} & {\small{}9} & {\small{}17} & {\small{}25} & {\small{}32}\tabularnewline
\hline 
\hline 
\multicolumn{13}{c}{{\small{}nearly symmetric systems}}\tabularnewline
\hline 
{\small{}rma10} & {\small{}4.85} & \textbf{\small{}2.31} & {\small{}2.53} & {\small{}5.02} & {\small{}2.47} & {\small{}2.69} & {\small{}67} & {\small{}9} & {\small{}9} & {\small{}3.4} & {\small{}2.0} & {\small{}2.3}\tabularnewline
\hline 
{\small{}PR02R} & {\small{}$-$} & \textbf{\small{}256} & {\small{}293} & {\small{}$-$} & {\small{}261} & {\small{}300} & {\small{}$-$} & {\small{}14} & {\small{}15} & {\small{}$-$} & {\small{}28} & {\small{}32}\tabularnewline
\hline 
\hline 
\multicolumn{13}{c}{{\small{}symmetric, saddle-point problems}}\tabularnewline
\hline 
{\small{}M3D3} & {\small{}$-$} & {\small{}53.6} & \textbf{52.9} & {\small{}$-$} & {\small{}77.3} & 76.8 & {\small{}$-$} & {\small{}52} & {\small{}53} & {\small{}$-$} & {\small{}2.6} & {\small{}2.6}\tabularnewline
\hline 
{\small{}M3D2} & {\small{}8.06} & {\small{}6.37} & \textbf{6.25} & {\small{}16.1} & {\small{}7.69} & 7.75 & {\footnotesize{}120} & {\small{}23} & {\small{}26} & {\small{}4.0} & {\small{}2.6} & {\small{}2.6}\tabularnewline
\hline 
\end{tabular}{\small\par}
\end{table}

\section{Conclusions and Future Work\label{sec:Conclusions and Future Work}}

In this paper, we described an MLILU technique, called HILUCSI, which
is designed for saddle-point problems from PDEs. The key novelty of
HILUCSI is that it takes into account the near or partial symmetry
of the underlying problems, and it improves the simplicity, robustness,
and efficiency of MLILU. More specifically, HILUCSI applies static
and dynamic deferring for improving robustness while enjoying a simpler
implementation than pivoting. It applies symmetric preprocessing techniques
at the top level for nearly or partially symmetric systems but applies
unsymmetric preprocessing and factorization at coarser levels, which
improved the robustness for problems from systems of PDEs. Furthermore,
the scalability-oriented dropping significantly improved the efficiency
of MLILU for large-scale problems. We demonstrated the robustness
and efficiency of HILUCSI as a right-preconditioner of restarted GMRES
for symmetric and unsymmetric saddle-point problems from mixed Poisson,
Stokes, and Navier-Stokes equations. Our results showed that HILUCSI
is significantly more robust than single-level ILU, such as SuperLU.
It also outperforms other MLILU packages (specifically, ILUPACK) by
a factor of four to ten for medium to large problems.

In its current form, HILUCSI has several limitations. First, if the
memory is severely limited, there may be too many droppings or too
many levels, and the preconditioner may lose robustness and efficiency.
We plan to optimize HILUCSI further for limited-memory situations.
Second, for vector-valued PDEs, the matrices may exhibit block structures.
It may be worthwhile to explore such block structures to improve cache
performance, similar to that in \cite{lishao10} and \cite{gupta2010adaptive}.
Finally, the HILUCSI algorithm is sequential as presented in this
paper. We will report its parallelization in the future.

\section*{Acknowledgments}

Results were obtained using the Seawulf and LI-RED computer systems
at the Institute for Advanced Computational Science of Stony Brook
University, which were partially funded by the Empire State Development
grant NYS \#28451. We thank Dr. Matthias Bollh\"ofer for helpful
discussions on ILUPACK. We thank the anonymous reviewers for their
helpful comments in improving the presentation of the paper. This
study does not have any conflicts to disclose.

\bibliographystyle{wileyj}
\bibliography{refs/refs,refs/multigrid,refs/compkrylov_refs,refs/psmilu,refs/refs_18}

\appendix

\section{\label{sec:Thresholding-in-inverse-based}Thresholding in inverse-based
dropping}

We motivate our thresholding strategy in inverse-based dropping within
each level of HILUCSI, based on a heuristic stability analysis to
bound $\rho\left(\boldsymbol{A}\boldsymbol{M}^{-1}-\boldsymbol{I}\right)$
(i.e., using the spectral radius as a ``pseudo-norm'' of $\boldsymbol{A}\boldsymbol{M}^{-1}-\boldsymbol{I}$).
Let $\boldsymbol{\delta}_{A}=\boldsymbol{A}-\boldsymbol{M}$, where
$\boldsymbol{M}=\boldsymbol{L}\boldsymbol{D}\boldsymbol{U}$, and
$\boldsymbol{A}=\left(\boldsymbol{L}+\boldsymbol{\delta}_{L}\right)\left(\boldsymbol{D}+\boldsymbol{\delta}_{D}\right)\left(\boldsymbol{U}+\boldsymbol{\delta}_{U}\right)$,
where $\boldsymbol{\delta}_{L}$, $\boldsymbol{\delta}_{D}$, and
$\boldsymbol{\delta}_{U}$ denote the perturbations to $\boldsymbol{L}$,
$\boldsymbol{D}$, and $\boldsymbol{U}$, respectively. Hence,
\begin{align*}
\boldsymbol{\delta}_{A} & =\left(\boldsymbol{L}+\boldsymbol{\delta}_{L}\right)\left(\boldsymbol{D}+\boldsymbol{\delta}_{D}\right)\left(\boldsymbol{U}+\boldsymbol{\delta}_{U}\right)-\boldsymbol{L}\boldsymbol{D}\boldsymbol{U}\\
 & =\boldsymbol{\delta}_{L}\boldsymbol{D}\boldsymbol{U}+\boldsymbol{L}\boldsymbol{\delta}_{D}\boldsymbol{U}+\boldsymbol{L}\boldsymbol{D}\boldsymbol{\delta}_{U}+\text{h.o.t.},
\end{align*}
where we omit the higher-order terms that involve more than one $\boldsymbol{\delta}$
matrix. Note that
\begin{align}
\rho\left(\boldsymbol{A}\boldsymbol{M}^{-1}-\boldsymbol{I}\right) & =\rho\left(\sqrt{\boldsymbol{D}^{-1}}\boldsymbol{L}^{-1}\boldsymbol{\delta}_{A}\boldsymbol{M}^{-1}\boldsymbol{L}\sqrt{\boldsymbol{D}}\right)\nonumber \\
 & \approx\rho\left(\sqrt{\boldsymbol{D}^{-1}}\left(\boldsymbol{L}^{-1}\boldsymbol{\delta}_{L}+\boldsymbol{\delta}_{D}\boldsymbol{D}^{-1}\right)\sqrt{\boldsymbol{D}}+\sqrt{D}\boldsymbol{\delta}_{U}\boldsymbol{U}^{-1}\sqrt{\boldsymbol{D}^{-1}}\right)\nonumber \\
 & \leq\sqrt{\kappa(\boldsymbol{D})}\left(\left\Vert \boldsymbol{L}^{-1}\boldsymbol{\delta}_{L}\right\Vert +\left\Vert \boldsymbol{\delta}_{D}\boldsymbol{D}^{-1}\right\Vert +\left\Vert \boldsymbol{\delta}_{U}\boldsymbol{U}^{-1}\right\Vert \right).\label{eq:modified-inverse-thresholding}
\end{align}
In dynamic deferring, we restrict the magnitude of the diagonal entries
to be no smaller than $1/\kappa_{D}$, and we estimate the maximum
magnitudes to be approximately equal to $\kappa_{D}$. Hence, $\sqrt{\kappa(\boldsymbol{D})}$
is bounded by $\kappa_{D}$, which leads to our thresholding for inverse-based
dropping in Section\ \ref{subsec:scalability-oriented-dropping}.
In terms of $\boldsymbol{L}^{-1}\boldsymbol{\delta}_{L}$ and $\boldsymbol{\delta}_{U}\boldsymbol{U}^{-1}$,
it is difficult to bound their 2-norms, and hence we approximately
bound $\left\Vert \boldsymbol{L}^{-1}\boldsymbol{\delta}_{L}\right\Vert _{\infty}$
and $\left\Vert \boldsymbol{\delta}_{U}\boldsymbol{U}^{-1}\right\Vert _{1}$
by $\left\Vert \boldsymbol{L}^{-1}\right\Vert _{\infty}\left\Vert \boldsymbol{\delta}_{L}\right\Vert _{\infty}$
and $\left\Vert \boldsymbol{U}^{-1}\right\Vert _{1}\left\Vert \boldsymbol{\delta}_{U}\right\Vert _{1}$,
respectively, as in \cite{bollhofer2001robust} and \cite{li2003crout}.
However, note that the thresholding strategy in \cite{bollhofer2001robust}
and \cite{li2003crout} did not take into account $\sqrt{\kappa(\boldsymbol{D})}$
(or $\kappa_{D}$). This omission was because they derived the thresholds
based on bounding $\boldsymbol{L}^{-1}\boldsymbol{A}\boldsymbol{U}^{-1}-\boldsymbol{D}$
instead of $\boldsymbol{A}\boldsymbol{M}^{-1}-\boldsymbol{I}$, probably
because it is impractical to bound $\boldsymbol{A}\boldsymbol{M}^{-1}-\boldsymbol{I}$
in a single-level ILU. In contrast, our derivation is for each level
of MLILU with dynamic deferring, so it is practical to bound some
norm of $\boldsymbol{A}\boldsymbol{M}^{-1}-\boldsymbol{I}$, which
corresponds to a stability measure of ILU \cite{benzi2002preconditioning}.

\section{\label{sec:Time-Complexity}Linear time complexity within each level}

For linear systems arising from PDEs, the input matrix $\boldsymbol{A}\in\mathbb{R}^{n\times n}$
typically has a constant number of nonzeros per row and per column.
We now show that the total cost of HILUCSI is linear in $n$ within
each level in this setting.

Let us analyze the cost of the first level in detail. Let $\boldsymbol{P}\in\mathbb{N}^{n\times n}$,
$\boldsymbol{Q}\in\mathbb{N}^{n\times n}$, $\boldsymbol{L}\in\mathbb{R}^{n\times m}$,
$\boldsymbol{D}\in\mathbb{R}^{m\times m}$ and $\boldsymbol{U}\in\mathbb{R}^{m\times n}$
be the output of the factorization of the current level. First, let
us consider the cost of updating $\boldsymbol{\ell}_{k}$ and $\boldsymbol{u}_{k}^{T}$
using the Crout version of ILU, or in short, \emph{Crout update}.
The total number of floating-point operations is bounded by
\begin{equation}
\mathcal{O}\left(\text{nnz}\left(\boldsymbol{L}+\boldsymbol{U}\right)\left(\max_{i\leq m}\left\{ \text{nnz}\left(\boldsymbol{a}_{i}^{T}\right)\right\} +\max_{j\leq m}\left\{ \text{nnz}\left(\boldsymbol{a}_{j}\right)\right\} \right)\right),\label{eq:bound-operations}
\end{equation}
where $\boldsymbol{a}_{i}^{T}$ and $\boldsymbol{a}_{j}$ denote the
$i$th row and $j$th column of $\boldsymbol{P}^{T}\boldsymbol{A}\boldsymbol{Q}$,
respectively.  Second, let us consider the cost of deferring and
dropping. Given an efficient data structure (see Section~\ref{subsec:Static-and-dynamic}),
the number of floating-point operations in dynamic deferring is proportional
to Crout update. Furthermore, in the scalability-oriented dropping,
we use quick select, which has an expected linear time complexity,
to find the largest nonzeros, followed by quick sort after dropping.
Hence, the time complexity of dropping is lower than that of Crout
update. If there is a constant number of nonzeros in each row and
column, then
\[
\max_{i\leq n}\left\{ \text{nnz}\left(\boldsymbol{a}_{i}^{T}\right)\right\} +\max_{j\leq n}\left\{ \text{nnz}\left(\boldsymbol{a}_{j}\right)\right\} =\mathcal{O}(1),
\]
and $\text{nnz}\left(\boldsymbol{L}+\boldsymbol{U}\right)=\text{nnz}\left(\boldsymbol{A}\right)=\mathcal{O}(n)$.
Hence, Crout update with deferring takes linear time. Third, for the
Schur component in \ref{eq:Schur-complement}, the most expensive
and also the most challenging part is the computation of $\boldsymbol{L}_{E}\boldsymbol{D}_{B}\boldsymbol{U}_{F}$.
We compute $\boldsymbol{L}_{E}$ and $\boldsymbol{U}_{F}$ along with
$\boldsymbol{L}_{B}$ and $\boldsymbol{U}_{B}$ during Crout update,
so its cost is bounded by (\ref{eq:bound-operations}). We compute
the sparse matrix-matrix multiplication (SpMM) using the algorithm
as described in \cite{bank1992sparse}. Since our scalability-oriented
dropping ensures that the nonzeros in each row of $\boldsymbol{L}_{E}$
and in each column of $\boldsymbol{U}_{F}$ are bounded by a constant
factor of those in the input matrix (see Section~\ref{subsec:scalability-oriented-dropping}),
the SpMM also takes linear time. As a side product, $\text{nnz}$
of $\boldsymbol{S}_{C}$ is linear in that of $\boldsymbol{A}$. 

For the other levels, the analysis described above also applies by
considering the fact that the nnz in the present level is proportional
to that in $\boldsymbol{A}$, and that the scalability-oriented dropping
is performed based on the nnz in each row and column in the input
matrix.

\end{document}